\documentclass[11pt,reqno,a4paper]{amsart}
\usepackage[T1]{fontenc}
\usepackage{lmodern,mathtools,amssymb,amsthm}
\usepackage[a4paper,left=26mm,right=26mm,top=26mm,bottom=28mm,headheight=12pt]{geometry}
\usepackage{booktabs,tabularx,microtype}
\usepackage[hidelinks]{hyperref}
\usepackage{textcase}
\hypersetup{pdftitle={A transversal Dirac theorem for Berge Hamilton cycles},pdfauthor={Ailian Chen, Nayi Jin, Zihao Wei},pdfsubject={A transversal Dirac theorem for Berge Hamilton cycles},pdfkeywords={Berge Hamilton cycle, r-graph system, transversal, minimum vertex degree},pdfdisplaydoctitle=true}
\newtheorem{theorem}{Theorem}[section]
\newtheorem{lemma}[theorem]{Lemma}
\newtheorem{proposition}[theorem]{Proposition}
\numberwithin{equation}{section}
\allowdisplaybreaks[1]
\newcommand{\ind}{\mathbf 1}
\DeclareMathOperator{\dist}{dist}

\title[Transversal Berge Hamilton cycles]{A transversal Dirac theorem for Berge Hamilton cycles}

\author[A. Chen, N. Jin, Z. Wei]{\NoCaseChange{Ailian Chen}\quad\textperiodcentered\quad
\NoCaseChange{Nayi Jin}\quad\textperiodcentered\quad
\NoCaseChange{Zihao Wei}\\[6pt]
\normalfont\small
\NoCaseChange{Department of Mathematics and Statistics, Fuzhou University, Fuzhou 350108, Fujian, China}\\[2pt]
\NoCaseChange{Corresponding author: Ailian Chen (\texttt{elian1425@fzu.edu.cn})}\\[2pt]
\NoCaseChange{Nayi Jin: \texttt{1262698242@qq.com}\quad
Zihao Wei: \texttt{2193318775@qq.com}}}
\date{}
\begin{document}

\begin{abstract}
Let $\mathcal H=(H_1,\ldots,H_n)$ be an $r$-graph system on $n$ vertices, $3\le r<n$. We prove that $\mathcal H$
has a rainbow Berge Hamilton cycle if every $H_i$ has minimum
vertex degree at least $\binom{\lfloor(n-1)/2\rfloor}{r-1}+1$ for
$2r<n$, or at least $r$ for $2r\ge n$. The bound is sharp for every
$(n,r)$, even for identical layers.
\end{abstract}

\maketitle

\begingroup
\small
\noindent\textit{2020 Mathematics Subject Classification.}\enspace
05C65, 05C45, 05C15.\par
\noindent\textit{Keywords.}\enspace
Berge Hamilton cycle, $r$-graph system, transversal, minimum vertex degree.
\endgroup
\medskip

\section{Introduction}
Dirac~\cite{Dirac} proved that every graph of order $n\ge3$ and
minimum degree at least $n/2$ has a Hamilton cycle. Aharoni asked whether the same bound in each of $n$ graphs on a common
$n$-vertex set guarantees a rainbow Hamilton cycle. Cheng, Wang and Zhao~\cite{CWZ}
proved an asymptotic version, and Joos and Kim~\cite{JK} proved the
exact result. Cheng et al.~\cite{CSWW} classified, for sufficiently large $n$,
systems with minimum degree at least $n/2-1$ in each member and
no rainbow Hamilton cycle. For asymptotic transfer results and other spanning subgraphs,
see \cite{GHMPS,MMP}.

Write $[n]=\{1,\ldots,n\}$. An \emph{$r$-graph} on $V$ is a family
of distinct $r$-element subsets of $V$. An \emph{$r$-graph system}
is an indexed family $\mathcal H=(H_1,\ldots,H_n)$ on a common
$n$-vertex set $V$. Its layers need not be distinct; the index $i$ is the colour of $H_i$. For an
$r$-graph $H$, put $d_H(v)=|\{h\in H:v\in h\}|$ and
$\delta_1(H)=\min_{v\in V}d_H(v)$.

In an $r$-graph $H$, a \emph{Berge path} of order $m\ge2$ and
length $m-1$ is a sequence
\[
 P=v_0h_0v_1h_1\cdots h_{m-2}v_{m-1},
\]
where $v_0,\ldots,v_{m-1}\in V(H)$ are distinct and the
representative hyperedges $h_0,\ldots,h_{m-2}\in H$ are distinct, with
$\{v_j,v_{j+1}\}\subseteq h_j$ for $0\le j\le m-2$.
A single vertex of $H$ is a Berge path of order one and length zero.
A \emph{Berge cycle} in $H$ of length $m\ge3$ is a sequence
\[
 C=v_0h_0v_1h_1\cdots v_{m-1}h_{m-1}v_0
\]
with distinct vertices in $V(H)$ and distinct representatives in
$H$ satisfying $\{v_j,v_{j+1}\}\subseteq h_j$ for every $j$
modulo $m$. Following~\cite{KLM}, let $V(Q)$ and $E(Q)$ denote
the core vertex set and representative hyperedge set of a Berge path
or cycle $Q$. Thus, if $P=v_0h_0\cdots h_{m-2}v_{m-1}$, then
$V(P)=\{v_0,\ldots,v_{m-1}\}$ and $E(P)=\{h_0,\ldots,h_{m-2}\}$;
if $C=v_0h_0\cdots v_{m-1}h_{m-1}v_0$, then
$V(C)=\{v_0,\ldots,v_{m-1}\}$ and
$E(C)=\{h_0,\ldots,h_{m-1}\}$. A representative in $E(Q)$ may
contain vertices outside $V(Q)$.
A Berge path or cycle $Q$ in $H$ is \emph{Hamilton} if $V(Q)=V(H)$.

For an $r$-graph system $\mathcal H=(H_1,\ldots,H_n)$ on $V$,
a Berge path or cycle in the $r$-graph $\bigcup_{i=1}^n H_i$ on $V$
is \emph{rainbow in $\mathcal H$} if its representatives admit
pairwise distinct colours $c_j\in[n]$ with $h_j\in H_{c_j}$.

Bermond, Germa, Heydemann and Sotteau~\cite{BGHS} studied degree
conditions for long Berge cycles. Coulson and
Perarnau~\cite[Corollary~1.5]{CP} obtained the sharp minimum vertex degree
threshold for Berge Hamilton cycles when $r=o(\sqrt n)$.
Related results concern degree conditions~\cite{MHG}, degree
sequences~\cite{Salia} and non-uniform hypergraphs~\cite{FKL}. Kostochka, Luo and
McCourt~\cite{KLM} determined the threshold for a single
$r$-graph for all $3\le r<n$. Throughout, put $k=r-1$,
$t=\lfloor(n-1)/2\rfloor$ and $\varepsilon=n-1-2t\in\{0,1\}$.
The threshold is
\begin{equation}\label{eq:threshold}
 D(n,r)=
 \begin{cases}
 \binom{t}{r-1}+1,&2r<n,\\
 r,&2r\ge n.
 \end{cases}
\end{equation}

The same threshold suffices in every layer of an $r$-graph system.

\begin{theorem}\label{thm:main}
Let $3\le r<n$, and let $\mathcal H=(H_1,\ldots,H_n)$ be an
$r$-graph system on $V$, where $|V|=n$.
If $\delta_1(H_i)\ge D(n,r)$ for every $i\in[n]$, then
$\mathcal H$ has a rainbow Berge Hamilton cycle. The bound is sharp for
every pair $(n,r)$.
\end{theorem}

For sharpness, take identical layers, as in~\cite{KLM}.
If $2r<n$, use two complete $r$-graphs on $\lceil n/2\rceil$
vertices, disjoint for even $n$ and sharing one vertex for odd $n$.
Their union has minimum degree $D(n,r)-1$ and a shadow that is
disconnected or has a cut vertex. If $2r\ge n$, delete one of
the $n$ cyclic intervals of length $r$. The remaining $n-1$ edges
have minimum degree $r-1$. Neither construction has a Berge
Hamilton cycle.

Bailey, Li and Luo~\cite{BLL}, together with the large-uniformity
extension of Bailey et al.~\cite{BHLL}, showed that the same
threshold guarantees Berge pancyclicity in a single $r$-graph for
sufficiently large $n$.
For tight Hamilton cycles, Cheng et al.~\cite{CHWWY} gave a transversal
codegree condition and Tang et al.~\cite{TWWY} a minimum $(r-2)$-degree
condition. For $3$-graph systems, the latter requires minimum vertex
degree at least $(5/9+\gamma)\binom n2$ in every layer, for fixed
$\gamma>0$ and sufficiently large $n$, whereas our Berge threshold
is $D(n,3)\sim\tfrac14\binom n2$.
Related work concerns stability~\cite{WGZ}, resilience in random
$r$-graphs~\cite{CEP} and rainbow Berge Hamilton cycles in randomly
edge-coloured random $r$-graphs~\cite{ZC}.

Coulson and Perarnau~\cite{CP} coloured each shadow edge by a
containing $r$-edge and used their rainbow Dirac theorem to obtain
distinct representatives along a Hamilton cycle. In a system,
different colours may require the same representative: $uv$ and
$vw$ may both require $uvw$. Thus colours and representatives
must be controlled separately.

We use the rotation framework and extremal choices of~\cite{KLM}.
In Section~\ref{sec:binomial}, Hall's theorem lifts ordinary rainbow
cycles with prescribed representatives (Lemma~\ref{merge:lifting});
reserving two representatives for a possible weak pair allows
rotations to preserve distinctness (Proposition~\ref{prop:even-one}).
For the remaining ranks, released neighbourhoods control exchanges
on Berge $(n-1)$-cycles. The final proof of Theorem~\ref{thm:main}
assembles these arguments with the long-cycle bound of
Section~\ref{sec:long}.

\subsection{General notation}\label{sec:notation}

We use $P_m$ and $C_s$ for ordinary paths and cycles on $m\ge1$
and $s\ge3$ vertices, respectively. From Section~\ref{sec:tools} onwards, Berge paths and cycles in
$\mathcal H$ are understood to be rainbow unless stated otherwise.
For distinct vertices $u,v$, we write $\{u,v\}$ for the $2$-set
when it is viewed as a subset of a hyperedge, and $uv$ for the
corresponding ordinary edge when it is viewed in a shadow or core graph.
The consecutive pairs of a Berge path or cycle form its ordinary
\emph{core graph}. For these Berge objects, adjacency, distance and
independence refer to the core graph.

On an oriented core cycle $C$, let $v^+,v^-$ be the successor and
predecessor of $v$, let $uCv$ denote the forward segment, and let
$C^{-1}$ denote the reverse orientation. Put
$A^\pm=\{v^\pm:v\in A\}$ and $N_C(A)=A^+\cup A^-$.
The trace of an edge $f$ on $C$ is $f\cap V(C)$. These conventions
also apply to oriented ordinary cycles.

For a colour $i$, put
\[
 H_i(v)=\{f\in H_i:v\in f\},\qquad
 S_i(\{u,v\})=\{f\in H_i:\{u,v\}\subseteq f\}.
\]
The shadow graph $F_i$ has edge set $\{uv:S_i(\{u,v\})\ne\varnothing\}$,
and $|S_i(\{u,v\})|$ is the codegree of $uv$ in $H_i$. Let
$\mathcal F=(F_1,\ldots,F_n)$. We write $uv:i$ when the ordinary
pair $uv$ is assigned colour $i$ in $\mathcal F$.
For selected representatives $R$ with assigned colours, write
\[
 d_R(v)=|\{h\in R:v\in h\}|,\qquad
 N_i^R(v)=\{u\ne v:\{u,v\}\subseteq f\text{ for some }f\in H_i\setminus R\},
\]
and let $\chi(R)$ be their assigned colour set. Write
$\overline{\chi(R)}=[n]\setminus\chi(R)$ for the set of unused colours.
The neighbourhood $N_i^R(v)$ excludes all of $R$, including
representatives assigned other colours. For an ordinary graph $G$,
let $c_j(G)$ be its number of $j$-cliques. Binomial coefficients
counting impossible choices are zero.

An ordinary path, cycle or forest is \emph{rainbow in $\mathcal F$}
if its pairs have distinct assigned colours and each colour-$i$ pair
belongs to $F_i$. A \emph{lifting} assigns pairwise distinct
representatives from the corresponding sets $S_i(\{u,v\})$; it
produces a Berge path, cycle or forest with the same core.

We record three elementary operations used repeatedly in the proofs. A
\emph{same-colour exchange} on a core pair $uv$ in colour $i$
replaces its current representative $g\in R\cap H_i$ by an unused
hyperedge $f\in H_i\setminus R$ containing $\{u,v\}$. The core pair
and its colour are unchanged, while $g$ becomes available. The change
in $d_R(z)$ is
$\mathbf1_{z\in f}-\mathbf1_{z\in g}$. In particular, if a cycle
maximizes $d_R(w)$ globally, an exchange whose old and new
representatives both avoid $w$ preserves this maximum and permits
comparison of a secondary extremum. Every later extremal comparison
requires all preceding extrema to be preserved; changing the omitted
vertex requires comparison with the global maximum.

For a Berge path $P=v_0h_0v_1\cdots h_{m-2}v_{m-1}$, an incidence
$v_0\in h_j$ permits a \emph{path rotation}: assign $h_j$ to the
core pair $v_0v_{j+1}$, retaining all other representatives and
colours. The new core order is
\begin{equation}\label{eq:berge-rotation}
 (v_j,v_{j-1},\ldots,v_0,v_{j+1},\ldots,v_{m-1}).
\end{equation}
Its representatives are $h_{j-1},\ldots,h_0,h_j,h_{j+1},\ldots,h_{m-2}$,
with their original colours.

\emph{Two-mark exclusion.} Let
$C=v_0h_0v_1h_1\cdots v_{s-1}h_{s-1}v_0$ be a rainbow Berge cycle,
let $w\notin V(C)$, and suppose no rainbow Berge cycle of length
$s+1$ exists. Mark $v_iv_{i+1}$ when $w\in h_i$, and orient the
cycle; the first and second vertices of an oriented marked pair are
its tail and head. Then every unused hyperedge in an unused colour
contains at most one marked tail and at most one marked head. Indeed,
if two marked pairs $v_iv_{i+1}$ and $v_jv_{j+1}$ occur in cyclic order
and an unused hyperedge contains $v_i,v_j$, that hyperedge, together
with the marked representatives on $v_{i+1}w$ and $wv_{j+1}$ and
the two arcs of $C$, gives the rainbow Berge cycle with core order
$v_i,v_j,\ldots,v_{i+1},w,v_{j+1},\ldots,v_i$.
This also applies when $v_{i+1}=v_j$. Reversing the orientation excludes
an unused hyperedge containing two marked heads.

\subsection{Cycles omitting one vertex}\label{sec:cycle-notation}

Let $\lambda$ be the maximum of $d_R(w)$ over rainbow Berge
$(n-1)$-cycles with representatives $R$ and omitted vertex $w$,
or zero if none exists. For a \emph{maximizing cycle} $C$, set
$R=E(C)$, $m=n-1$, and orient its core as $(v_0,\ldots,v_{m-1})$. Write $e_i=v_iv_{i+1}$ for a core
pair and $h_i\in H_{c_i}$ for its representative, with indices modulo
$m$. A core pair is \emph{marked} if $w\in h_i$; let $Y^+$ be the
set of marked heads. For the unused colour $\alpha$, put
\[
 \mathcal T=\{f\setminus\{w\}:f\in H_\alpha(w)\setminus R\},
 \qquad U=\bigcup_{T\in\mathcal T}T.
\]
Let $\Gamma$ be the graph on $U$ joining vertices in a common trace.
A vertex is \emph{active} if it belongs to $U$. For $v_i\in U$, the
neighbourhood released at $v_{i+1}$ in colour $c_i$ is
$T_i=N_{c_i}^{\,R\setminus\{h_i\}}(v_{i+1})$.

\section{Links and Berge path rotations}\label{sec:tools}

\subsection{Elementary counting bounds}

\begin{lemma}\label{lem:link-incidence}
Let $\mathcal L\subseteq\binom{[N]}q$, $q\ge2$, $q-1\le s_0<N$.
Let $\beta\ge0$. If at most $s_0$ vertices have degree greater than
$\beta$ in the $q$-graph $\mathcal L$, then
$|\mathcal L|\le\binom{s_0}q+\beta(N-s_0)$.
\end{lemma}
\begin{proof}
Choose an $s_0$-set $X$ containing every vertex of degree above $\beta$.
At most $\binom{s_0}q$ members lie in $X$; each other member
uses one of the at most $\beta(N-s_0)$ incidences outside $X$.
\end{proof}

\begin{lemma}\label{lem:separation}
For an ordinary cycle $C_s$, $s\ge3$, and an ordinary path $P_m$,
$m\ge1$, the following hold.
\begin{enumerate}
\item[(i)] If $A,B\subseteq V(C_s)$, $A\ne B$, $|A|=|B|=d$ and
$s\le2d+1$, then some distinct $a\in A$, $b\in B$ satisfy
$ab\in E(C_s)$.
\item[(ii)] Let $A,B\subseteq V(C_s)$ have $|A|,|B|\ge2$, and let
$\ell\ge2$. If $\dist_{C_s}(a,b)\ge\ell+1$ for every
$a\in A$, $b\in B$ with $a\ne b$, then
$s\ge|A|+|B|+2\ell$ when $A\ne B$, and
$s\ge(\ell+1)|A|$ when $A=B$.
\item[(iii)] If $A,B\subseteq V(P_m)$ and no edge of $P_m$ has one
end in $A$ and the other in $B$, then
$|A|+|B|\le m+\ind_{m\text{ odd}}$.
If $|A|,|B|\ge\lceil m/2\rceil$ as well, then $A=B$ is a maximum
independent set of $P_m$. When $m$ is odd this set is unique.
\item[(iv)] Let $t\ge6$, $3\le q\le t$, $1\le p\le q$, and
$t+2\le s\le2t+2-q$. Put $\ell=p-\mathbf1_{p=q}$. If
$A,B\subseteq V(C_s)$ satisfy $|A|,|B|\ge t-p+1$ and
$\dist_{C_s}(a,b)\ge\ell+1$ for every distinct
$a\in A,b\in B$, then $p=q=t$ and $s=t+2$.
\end{enumerate}
\end{lemma}
\begin{proof}
List $A\cup B$ in cyclic order, taking the $a=|A\cap B|$ common
vertices as singleton blocks and the maximal runs in $A\setminus B$
or $B\setminus A$ as the other $b$ blocks. If distinct cross
vertices have distance at least $\ell+1$ and there are at least
two blocks, each gap between blocks has at least $\ell$ vertices.
Thus
\begin{equation}\label{eq:cyclic-block-gaps}
 s\ge |A|+|B|+(\ell-1)a+\ell b.
\end{equation}
For (i), take $\ell=1$; distinct equal-sized sets have $b\ge2$,
contradicting $s\le2d+1$. For (ii), if neither set contains the other, then $b\ge2$. Under strict containment, $a\ge2$, $b\ge1$
and $(\ell-1)a+\ell b\ge2\ell$. If $A=B$, sum the gaps between
successive vertices of $A$.

For (iii), assume $A\cup B\ne\varnothing$; otherwise the claims
hold. The linear count gives $m\ge |A|+|B|+b-1$. If $A\ne B$, then $b\ge1$, giving $|A|+|B|\le m$.
If both sizes are at least $\lceil m/2\rceil$, this is impossible
for odd $m$. For even $m$, both sizes would be $m/2$, forcing both set
differences to be nonempty and $b\ge2$, again a contradiction.
Thus $A=B$, an independent set of size at most $\lceil m/2\rceil$.
On an odd path the maximum independent set consists uniquely of
alternate vertices starting at an endpoint.

For (iv), first suppose $p=1$. Then $\ell=1$ and
$s\le2t-1$. If $A\ne B$, choose distinct $t$-subsets
$A'\subseteq A$ and $B'\subseteq B$; part (i) gives an adjacent
cross-pair, a contradiction. If $A=B$, then $A$ is independent
and $|A|\ge t>\lfloor s/2\rfloor$, also impossible.

Suppose $2\le p<q$. Then $\ell=p$. If $A\ne B$, part (ii)
gives $s\ge2(t-p+1)+2p=2t+2$, contrary to
$s\le2t+2-q\le2t-1$. If $A=B$, part (ii) gives
$s\ge(p+1)(t-p+1)>2t+1-p\ge s$, since
$(p+1)(t-p+1)-(2t+1-p)=(p-1)(t-p)>0$.

Finally suppose $p=q<t$. Since $q\ge3$, we have $p\ge3$ and
$\ell=p-1$. If $A\ne B$, part (ii) gives $s\ge2t$, while
$s\le2t+2-p\le2t-1$. If $A=B$, part (ii) gives
$s\ge p(t-p+1)>2t+2-p\ge s$, because
$(p-2)(t-p)-2\ge t-5>0$. The only remaining case is
$p=q=t$, and the bounds on $s$ then give $s=t+2$.
\end{proof}

\begin{lemma}\label{lem:cyclic-gap}
Let $S$ be an independent set of at least $q\ge1$ vertices in an
ordinary cycle $C_s$. Every gap between consecutive vertices of $S$,
measured in core edges, has length at most $s-2(q-1)$.
\end{lemma}
\begin{proof}
The vertices of $S$ divide the cycle into at least $q$ gaps. Since
$S$ is independent, every gap other than a fixed one has length at
least two. The sum of all gap lengths is $s$, giving the stated bound.
\end{proof}

\subsection{Path rotations and endpoint neighbourhoods}

\begin{lemma}\label{lem:common-endpoint-cuts}
Let $\mathcal H=(H_1,\ldots,H_n)$ be an $r$-graph system on $n$ vertices,
$3\le r<n$. Let $P=v_0h_0v_1\cdots h_{m-2}v_{m-1}$ be a rainbow Berge path,
$m\ge3$, with representatives $R$.
For distinct unused colours $a,b\in\overline{\chi(R)}$, put
\[
 I=\{j:v_{j+1}\in N_a^R(v_0)\},\qquad
 J=\{j:v_j\in N_b^R(v_{m-1})\},\quad 0\le j<m-1.
\]
\begin{enumerate}
\item[(i)] If $\mathcal H$ has no rainbow Berge $m$-cycle, then
$v_{m-1}\notin N_a^R(v_0)$, $v_0\notin N_b^R(v_{m-1})$,
$I\cap J=\varnothing$ and
\[
 |N_a^R(v_0)\cap V(P)|+|N_b^R(v_{m-1})\cap V(P)|\le m-1.
\]
If also $m=n-1$ and no rainbow Berge Hamilton cycle exists,
the omitted vertex belongs to at most one of
$N_a^R(v_0),N_b^R(v_{m-1})$.
\item[(ii)] Suppose $V=V(P)\dot\cup\{y\}$ and no rainbow Berge
Hamilton cycle exists. If every rainbow Berge cycle on $V(P)$ has
at most $\lambda$ representatives containing $y$, but $P$ has at
least $\lambda+1$, then
\[
 |N_a^R(v_0)|+|N_b^R(v_{m-1})|+\kappa\le m,
\]
where $\kappa$ counts the indices outside $I\cup J$.
\end{enumerate}
\end{lemma}
\begin{proof}
An unused hyperedge of either colour containing $\{v_0,v_{m-1}\}$
would close $P$. This contradicts the forbidden $m$-cycle in (i),
or the bound on the number of representatives containing $y$ in
(ii). Hence
$v_{m-1}\notin N_a^R(v_0)$ and $v_0\notin N_b^R(v_{m-1})$.
Moreover, two attachments in colours $a,b$ at opposite ends are
distinct: a common physical representative would contain both
$v_0$ and $v_{m-1}$.

For (i), if $j\in I\cap J$, choose unused hyperedges in colours
$a,b$ containing $\{v_0,v_{j+1}\}$ and
$\{v_j,v_{m-1}\}$, respectively. Together with the path outside
$h_j$, they form a rainbow Berge $m$-cycle, a contradiction.
Thus $I\cap J=\varnothing$. The endpoint exclusions show that
$N_a^R(v_0)\cap V(P)=\{v_{j+1}:j\in I\}$ and
$N_b^R(v_{m-1})\cap V(P)=\{v_j:j\in J\}$, so their total size
is $|I|+|J|=|I\cup J|\le m-1$.
If $m=n-1$, the only vertex outside $P$ is $y$. If it belonged
to both neighbourhoods, two unused attachments would insert $y$
between the two ends of $P$, giving a rainbow Berge Hamilton cycle.
Hence it belongs to at most one neighbourhood.

For (ii), every common cut $j\in I\cap J$ must be marked, that is,
$y\in h_j$. Otherwise replace $h_j$ by the two unused attachments
on $v_0v_{j+1}$ and $v_jv_{m-1}$. The resulting cycle on $V(P)$
retains all $\lambda+1$ representatives containing $y$, a
contradiction.

There cannot be two common cuts. Indeed, if $j<k$ both lie in
$I\cap J$, use $h_j,h_k$ on $v_jy,yv_{k+1}$, and use the unused
$b$-edge on $v_{m-1}v_k$ and the unused $a$-edge on
$v_{j+1}v_0$. The resulting Hamilton core order is
\[
 v_0\,P\,v_j,\ y,\ v_{k+1}\,P\,v_{m-1},\
 v_k\,\overleftarrow P\,v_{j+1},\ v_0,
\]
where a segment of length zero is omitted. If there is exactly
one common cut $j$ and $y\in N_b^R(v_{m-1})$, use $h_j$ on
$v_jy$, an unused $b$-edge on $yv_{m-1}$, and an unused $a$-edge
on $v_{j+1}v_0$ to obtain a Hamilton cycle. The reverse argument
applies if $y\in N_a^R(v_0)$. Finally, if $y$ belongs to both
neighbourhoods, the two distinct unused attachments on
$v_0y$ and $yv_{m-1}$ close $P$ to a Hamilton cycle.

Put
$e=\mathbf1_{y\in N_a^R(v_0)}+\mathbf1_{y\in N_b^R(v_{m-1})}$.
The preceding cases give $|I\cap J|\le1-e$, and exactly
$\kappa$ indices lie outside $I\cup J$. Thus
\[
 |N_a^R(v_0)|+|N_b^R(v_{m-1})|
 =|I|+|J|+e
 \le m-1-\kappa+(1-e)+e=m-\kappa.
\]
This proves (ii).
\end{proof}

\begin{lemma}\label{lem:fixed-end-head-bound}
Let $\mathcal H=(H_1,\ldots,H_n)$ be an $r$-graph system on $n$ vertices,
$3\le r<n$. Let $P$ be a rainbow Berge path with core vertex set
$V(P)$, $|V(P)|=m\ge3$, ends $x,z$ and representatives $R$, and let
$y\in V\setminus V(P)$. Suppose every rainbow Berge cycle on $V(P)$ has
at most $\lambda$ representatives containing $y$, but $P$ has at
least $\lambda+1$. For every unused colour $a$,
$|N_a^R(x)\cap V(P)|+d_R(z)\le m-1$.
\end{lemma}
\begin{proof}
Rotations fixing $x$ preserve $R$ and its colours and make the
$d_R(z)$ heads of representatives containing $z$ attainable ends.
None lies in $N_a^R(x)$, since closure would retain at least
$\lambda+1$ representatives containing $y$. Both sets lie in
$V(P)\setminus\{x\}$, proving the bound.
\end{proof}

\subsection{A long-cycle bound}

We use a rainbow version of the circumference bound in~\cite{KLM}.

\begin{lemma}\label{lem:all-circumference}
Let $\mathcal H=(H_1,\ldots,H_n)$ be an $r$-graph system on $n$ vertices.
Let $t$ be an integer with $3\le r\le t\le n-2$.
If $\delta_1(H_i)>\binom t{r-1}$ for every $i\in[n]$, then
$\mathcal H$ has a rainbow Berge cycle of
length at least $t+2$.
\end{lemma}
\begin{proof}
Suppose every rainbow Berge cycle has length at most $t+1$.
A longest Berge path has at least $t+1$ vertices. Indeed, if it
had at most $t$, then in an unused colour its end would be incident
with at most $\binom{t-1}{r-1}$ unused incident hyperedges and at
most $t-1$ selected representatives. Its degree would then be at
most $\binom{t-1}{r-1}+t-1\le\binom t{r-1}$, a contradiction.
Among oriented longest paths $P=v_0h_0v_1h_1\cdots$, maximize the
number of $h_0,\ldots,h_{t-1}$ contained in $X=\{v_0,\ldots,v_t\}$.
Put $R=E(P)$. Every unused edge containing $v_0$ in an unused
colour lies in $X$; otherwise it extends $P$ or closes an initial
segment of length at least $t+2$. The same argument excludes
$v_0$ from all representatives after $h_t$.

Fix an unused colour $a$. If no edge of $H_a\setminus R$ contains
$v_0$, then $r=t$ and all $h_0,\ldots,h_t$ contain $v_0$.
Rotations at these edges make every vertex of $X$ initial while
preserving $X$, the first $t$ representatives and the secondary maximum. If none of these vertices had an unused colour-$a$ edge,
the degree bound would force $X\subseteq h_0$, impossible.
Choose this orientation and put
$J=\{j\in[t]:v_j\in N_a^R(v_0)\}$.
An unused edge on $v_0v_j$ makes $v_{j-1}$ the initial vertex and
releases $h_{j-1}$, which lies in $X$ by maximality.
For $r=t$, one edge gives $|J|\ge t-1$. For $r=3<t$,
there are at least $\binom t2-t$ unused colour-$a$ hyperedges
through $v_0$, more than $\binom{t-2}2$; hence $|J|\ge t-1$.
For $4\le r<t$, suppose $j_0\in[t]\setminus J$. Then every
unused colour-$a$ hyperedge through $v_0$ is contained in
$X\setminus\{v_{j_0}\}$, so there are at most
$\binom{t-1}{r-1}$ such hyperedges, but at least
$\binom t{r-1}-t$ unused hyperedges through $v_0$.
Since $\binom{t-1}{r-2}\ge\binom{t-1}{2}>t$, we have
$\binom t{r-1}-t>\binom{t-1}{r-1}$, a contradiction.
Thus $J=[t]$.

Suppose $r=t$ and some $h_e$, $e<t$, leaves $X$. Its head
$p=v_{e+1}$ is not an unused neighbour of $v_0$ in any unused
colour, since rotation there would improve the maximum.
The family
\[
 \{g\in H_a\setminus R:v_0\in g\}
\]
is nonempty by the preceding choice of orientation. Every member
lies in $X$ and avoids $p$. Since $|X|=t+1=r+1$, the only possible
member is $f=X\setminus\{p\}$; hence the family is exactly
$\{f\}$.
For every $j\in[t]\setminus\{e+1\}$, rotation using $f$ on $v_0v_j$
replaces an edge inside $X$ and preserves the maximum. The released
$h_{j-1}$ lies in $X$ and avoids the new head of $h_e$, either $p$
or $v_e$. Since $h_{j-1}\ne f$, it equals $X\setminus\{v_e\}$.
This makes $t-1\ge2$ representatives identical, a contradiction.

For each $j\in J$, the rotation to $v_{j-1}$ preserves $X$ and the
secondary maximum. Every unused incident edge in the released colour
lies in $X$, so degree requires a selected edge at $v_{j-1}$ leaving
$X$, among the first $t+1$ representatives. The boundary edge $h_t$, if present, contains $v_t,v_{t+1}$,
so covers at most $r-2<|J|$ of these initial vertices.
Only when $r=3<t$ can one of the first $t$ representatives $h_e$
also leave $X$. It contains at most one of these initial vertices:
$v_e$ is excluded because $e+1\notin J$, and another vertex lies
outside $X$. The two triples cover at most two of $|J|\ge t-1\ge3$ initial
vertices. In each case an initial vertex has no selected edge
leaving $X$, contradicting degree.
\end{proof}

\section{Berge cycles of length at least \texorpdfstring{$n-1$}{n-1}}\label{sec:long}
\begin{lemma}\label{r2:triple-two-omitted}
Let $\mathcal H=(H_1,\ldots,H_n)$ be an $r$-graph system on
$n=2t+1+\varepsilon$ vertices, $\varepsilon\in\{0,1\}$, with
$3\le r\le t$ and $\delta_1(H_i)\ge D(n,r)$ for every $i$.
If $\mathcal H$ has a rainbow Berge $(n-2)$-cycle, then it has a
rainbow Berge cycle of length at least $n-1$.
\end{lemma}
\begin{proof}
Suppose no rainbow Berge cycle has length at least $n-1$.
Choose an $(n-2)$-cycle with representatives $R$ and ordered omitted
vertices $u,v$ maximizing $(d_R(u)+d_R(v),d_R(u))$ lexicographically,
subject to $d_R(u)\ge d_R(v)$. Put $s=n-2$.

Write the core cycle as
$C=v_0h_0v_1h_1\cdots v_{s-1}h_{s-1}v_0$, with indices modulo $s$.
For either omitted vertex $w\in\{u,v\}$, mark $v_iv_{i+1}$ when
$w\in h_i$. The two-mark exclusion applies to these marked pairs,
since no rainbow Berge cycle of length $s+1=n-1$ exists.

\emph{Restrictions on traces.} For each unused colour, the trace
on the core cycle of an unused edge containing $u$ is independent
and avoids the ends of all core edges marked at $u$. If such an edge $f$ contains a core pair represented by $g$, insert
$u$ when $u\in g$, or $v$ when $v\in f\cap g$. Otherwise replace
$g$ by $f$ in its unused colour, increasing $(d_R(u)+d_R(v),d_R(u))$
lexicographically. Meeting an end of a $u$-marked pair likewise permits insertion of $u$.
Consequently, no core edge joins the trace union in one unused
colour to the trace union in the other. The representatives are distinct, since a common one would
contain a core edge in its trace. If $d_R(u)=0$, then $d_R(v)=0$, and the degree bound forces both trace
unions to have size at least $t$. After cutting a core pair, Lemma~\ref{lem:separation}(iii)
makes the two unions coincide. Their common set $\Phi$ is independent
in $C$ by the adjacency exclusion, so $s=2t$ and $\Phi$ alternates. The degree bound also gives an unused edge $f$ containing $uv$,
since edges containing $u$ and avoiding $v$ have at most
$\binom t{r-1}$ choices. As $0<|f\cap \Phi|=r-2<t$, there are
successive vertices $x,y$ in the cyclic order on $\Phi$ with $x\in f$
and $y\notin f$. Replace the core vertex between $x$ and $y$ by
$u$, using $f$ on $xu$ and an unused edge of the other colour on
$uy$. These representatives are distinct because $y\notin f$.
The new cycle still omits $v$ but uses $f$, with $v\in f$,
increasing $d_R(u)+d_R(v)$.
Hence $d_R(u)>0$.

Suppose each unused colour has an unused edge containing
$uv$. Then the traces of unused edges containing $v$ are
independent, avoid the ends of all core edges marked at either
omitted vertex, and have no adjacent vertices in different
colour trace unions. If $xy\subseteq f$ and $v\in f$, insert $v$ if the core pair's
representative contains $v$; if it avoids $u,v$, exchange it for $f$
to increase $d_R(u)+d_R(v)$. If $xy$ is marked at $u$, it is enough that $y\in f$: the
core edge $xy$ can be replaced by $x,u,v,y$, using its marked
representative, an unused edge of the other colour containing
$uv$, and $f$. The two unused representatives differ because $u\notin f$.

\emph{Unused edges at the omitted vertices.} Suppose $d_R(u)\ge D(n,r)$. Then $r=t$ or $(r,t)=(3,4)$,
since otherwise $D(n,r)>s$. In either case $s-d_R(u)\le r-1$. Each unused colour
has at least $\lceil nD(n,r)/r\rceil-s\ge5$ unused edges. By the two-mark exclusion, each such edge contains at most one
marked tail and at most one marked head. An unused edge containing $u$ needs $r-2$ vertices outside the
marked ends. If every core pair were marked, no such edge could
exist; an unused edge avoiding $u$ would contain two marked tails.
Since unused edges exist, the marks do not cover the cycle. Let $c$ count the maximal paths of marked core edges. The set $O$ of core vertices outside their
ends has size $s-d_R(u)-c\le r-1-c$. If $c\ge2$, an unused edge avoids $u$ and needs at least $r-1$
core vertices, chosen from $O$ together with at most one marked
tail and at most one marked head. Thus $c=2$, and the edge must
consist of $v$, all vertices of $O$, one end that is only a
marked tail, and one end that is only a marked head. Each endpoint has at most two choices, giving at most four unused
edges, contrary to the lower bound five.

At $c=1$ the unmarked segment of the ordinary core cycle has at most $r$
vertices. An unused edge containing $v$ and avoiding $u$ cannot contain an
unmarked core pair, since insertion or a representative exchange
would contradict the extremal choice. If an unused edge
contains $u$ or an internal vertex of the marked path, it must
contain $v$ and all $r-2$ vertices of $O$: an internal vertex is
both a marked head and a marked tail. For $r\ge4$ this covers an unmarked pair.
The remaining edges lie in the unmarked core segment, or meet it in an
independent $(r-1)$-set and contain $v$. The latter is impossible
for $r\ge4$, and even for triples there are at most two choices in total.
Thus $r=3$, $|O|=1$ and $(n,d_R(u))=(8,4)$. Each of the three internal
marked vertices lies in at most one unused triple, with $v$ and $O$.
Degree requires nine selected incidences there, but the four marks
contribute six and the other two selected triples at most one each. Hence $d_R(u)<D(n,r)$. Since $d_R(v)\le d_R(u)$, each unused colour
has an unused edge at each omitted vertex.

\smallskip\noindent\emph{The sizes of the trace unions.}
Let $A,B$ be the trace unions at $u$ in the two unused colours,
labelled so that $d=|A|\le|B|$. The oriented core edges
leaving $A$ and entering $B$ form disjoint sets, and none is marked at $u$. Hence
$k-1\le d$ and $D(n,r)\le d_R(u)+\binom{d+1}k\le s-2d+\binom{d+1}k$.
On the integer interval $k-1\le d\le t-2$, the last expression
is convex and has endpoint values
$2t+2+\varepsilon-2k$ and $3+\varepsilon+\binom{t-1}k$.
Both values are less than $D(n,r)$ when $t\ge5$, since $\binom{t-1}{k-1}\ge t-1$,
and $\binom tk\ge\binom t2$ unless $k=t-1$, when $D(n,r)=t+1$.
For $t=3,4$ the bounds leave $r=3,d=1$ or $(t,\varepsilon,d)=(4,1,2)$.
In the former $d_R(u)=D(n,r)-1\ge s-3$, so at most two vertices avoid the
marked ends, adjacent if there are two. In the latter $d_R(u)=4$, every
$r$-edge containing $u$ and contained in $\{u,v\}\cup A$ occurs, and the four pairs avoiding
the independent set $A$ are marked. In both cases the vertices outside
the marked ends lie in $A\cup N_C(A)$. Thus $B=A$, and each unused colour contains all $\binom{d+1}k=D(n,r)-d_R(u)$ possible
$r$-edges containing $u$ and contained in $\{u,v\}\cup A$. The edges containing $uv$ put $A$ in both trace unions at $v$;
the adjacency and marked-end exclusions force equality. Degree gives $D(n,r)\le d_R(v)+\binom{d+1}k=d_R(v)+D(n,r)-d_R(u)$,
so $d_R(v)=d_R(u)$.
For triples, the two sets of marked core edges are disjoint,
since a triple cannot contain both omitted vertices and a core
pair. This contradicts $2d_R(u)>s-2d$. For $r=4$, the same four pairs avoiding $A$ are marked at both
omitted vertices. Two are consecutive on the six vertices outside
$A$: replace their segment $a,b,c$ by $a,u,v,c$, using the two
marks and an unused edge containing $uv$. Thus $d\ge t-1$.

\smallskip\noindent\emph{Replacing a core vertex.}
Delete the ends of a $u$-marked core pair from the ordinary core
cycle. On the remaining ordinary path $P_{2t-3+\varepsilon}$,
Lemma~\ref{lem:separation}(iii) gives $A=B=S$, where $S$ is
an independent set of size $t-1$. Hence $d_R(u)\le s-2|S|=1+\varepsilon$, and in each unused colour at most $d_R(u)-1\le1$ of the $r$-edges
containing $u$ and contained in $\{u,v\}\cup S$ are absent.
For $r=3$, we claim that $d_R(u)=1$. Each unused colour has an
unused edge containing $uv$, because $|S|\ge2$ and at most one
of the possible triples is absent. The trace exclusions at $v$,
together with $D(n,r)-d_R(v)>\binom{t-1}2$, imply that both trace unions
at $v$ have size at least $t-1$. Deleting the ends of a core edge
marked at $u$ and applying Lemma~\ref{lem:separation}(iii) again
gives a common independent set $T$ of size $t-1$. The degree
bound forces $d_R(v)\ge1$. The $d_R(u)+d_R(v)$ marked core edges are distinct
and avoid $T$, so $d_R(u)+d_R(v)\le s-2(t-1)\le2$. Since $d_R(u)\ge d_R(v)$,
we obtain $d_R(u)=d_R(v)=1$.

At most two core edges avoid $S$. Since at least one is marked
at $u$, some marked edge $ab$ is adjacent to a vertex $x\in S$;
reverse the orientation if necessary so that the order is $x,a,b$.
Choose distinct unused edges in the two unused colours, containing
$x,u,v$ and $u,v$, respectively. For triples, every edge
containing $u$ and contained in $\{u,v\}\cup S$ is available
in each unused colour, and $|S|\ge2$ ensures distinct choices.
For $r\ge4$, the candidate families have sizes
$\binom{t-2}{r-3}\ge2$ and $\binom{t-1}{r-2}\ge3$,
respectively, with at most one missing member each; distinct
choices again exist. Replace $x,a,b$ by $x,v,u,b$, using
the old representative of $ab$ on $ub$, to obtain a Berge $(s+1)$-cycle.
\end{proof}

\begin{proposition}\label{prop:boundary-long}
Let $\mathcal H=(H_1,\ldots,H_n)$ be an $r$-graph system on
$n=2t+1+\varepsilon$ vertices, $\varepsilon\in\{0,1\}$, with
$3\le r\le t$ and $\delta_1(H_i)\ge D(n,r)$ for every $i$.
Then $\mathcal H$ has a rainbow Berge cycle of length at least $n-1$.
\end{proposition}
\begin{proof}
Suppose the maximum rainbow Berge cycle length is $s<n-1$.
Lemmas~\ref{lem:all-circumference} and~\ref{r2:triple-two-omitted}
give $t+2\le s\le n-3$. Put $q=n-s$ and $q_0=q-\ind_{r<t}$.
Choose a Berge $s$-cycle $C$ and a Berge path $P$ of order at most
$q_0$, disjoint in core vertices, representatives and colours,
maximizing first $h=|V(P)|$, then $\sum_{f\in E(C)}|f\cap V(P)|$.
The path may have order one, and its representatives may meet $V(C)$.
Set $R=E(C)\cup E(P)$; $q-h+1$ colours are unused. Orient $C$.

Let $x,z$ be the ends of $P$ and $c(w)=|\{f\in E(C):w\in f\}|$.
Put $W=V(P)$ if $h<q_0$, and $W=V\setminus V(C)$ otherwise;
write $p=|W|$. Traces are taken on $C$.

If $h=1<q_0$, no unused edge in an unused colour contains two
vertices outside $C$, since they would extend $P$. In general,
every such edge $f$ with $x\in f$ lies in $W\cup V(C)$.
If its trace contains a core pair represented by $g$, either
a vertex $y\in g\cap V(P)$ allows a longer detour along $xPy$, with attachments
$f,g$, or replacing $g$ by $f$ improves the second extremum.
Thus the trace is independent. It also avoids the ends of every
$x$-marked pair, or $x$ could be inserted. The same holds at $z$.

Distinct unused attachments in different colours at $x,z$ with
distinct core endpoints require distance at least $h+1$, or the
detour along $P$ lengthens $C$. A mark at $x$ also supplies an
attachment, so each opposite trace vertex excludes $2h$ consecutive
mark indices. For a nonempty opposite trace union $B$, their union
has size at least $\min\{s,|B|+2h-1\}$, giving
$c(x)\le\max(0,s-|B|-2h+1)$.
Likewise, distinct tails marked at opposite ends require distance
at least $h$: the arc spanning their pairs has one more edge.

\emph{Suppose first that $r<t$.}
Here $h=p-\ind_{p=q}$. Fix unused colours $a\ne b$ at $x,z$. We show that each end has
at least $t-p+1$ core neighbours in its assigned colour, each with
two unused representatives when $r\ge4$.

First let $r=3$, and write
$A=N_a^R(x)\cap V(C)$ and $B=N_b^R(z)\cap V(C)$. If $p=1$
and $A=\varnothing$, then $D(n,r)\le s\le2t-1\le D(n,r)$, so every
cycle representative contains $x$. Degree summation gives an unused edge with two core vertices, hence
two marked tails, contradicting the two-mark exclusion. Thus
$A\ne\varnothing$, and each unused edge at $x$ has two core
vertices, so $d=|A|\ge2$. Let $e_C(A)$ denote the number of core
edges with both ends in $A$. At most $s-2d+e_C(A)$ core edges
avoid $A$, and only these can be marked at $x$. Trace independence
also excludes those $e_C(A)$ pairs from the unused link. Hence
$D(n,r)\le\binom d2+s-2d$. Convexity rules this out for
$2\le d\le t$, so $|A|\ge t+1$. The same argument applies to
$B$.

When $p\ge2$, every selected triple containing $x$, except
a representative of $C$ or the first representative of $P$, lies
in $W$: its other two vertices are the ends of an internal core
edge of $P$. Thus
$D(n,r)\le\binom{p-1+|A|}2+c(x)+1$. If both unions have size at most $t-p$,
the two disjoint sets of marked core edges each have size at
least $t-1$. They cover all but at most one core edge, so two
edges marked at opposite ends are consecutive, contrary to the
distance restriction on their tails. If only $|A|\le t-p$,
then $c(x)\ge t-1$, while $|B|\ge t-p+1$ gives by the arc bound
$c(x)\le t+2-p-q$ for $p<q$, or $c(x)\le t+4-2p$ for $p=q$.
Both upper bounds are less than $t-1$, a contradiction.

For $r\ge4$  and $K=\binom{t-1}{k-1}$. First assume
\[
 s+h-1+\binom{p-1}k<D(n,r),\qquad s+h+p<K+t+3.
\] The first inequality makes the trace union at $x$ nonempty in
each unused colour. Its size $d$ is less than $s$, since otherwise
a trace in the other unused colour would give adjacent
attachments that insert $x$. Let $j$ of its $d$ neighbours in colour $a$ have at least two unused
representatives. At least $d+1$ core pairs meet the union and are
unmarked. The other $d-j$ neighbours account for at most $d-j$
unused edges; the remaining edges at $x$ lie in $W$ together with
those $j$ neighbours. Hence
\[
 D(n,r)\le (s+h-d-2)+\binom{p-1+j}k+d-j.
\]
On $0\le j\le t-p$ the right side is convex and both endpoint values
are below $D(n,r)$ by the displayed conditions. Thus $j\ge t-p+1$.
Both inequalities hold for $t\ge7$, since
$s+h\le2t+1$, $p\le t$ and $K\ge\binom{t-1}2$. For $t=5,6$, substituting
$3\le q\le t-1+\varepsilon$, $s=2t+1+\varepsilon-q$ and
$p\in\{1,\ldots,q-2,q\}$ leaves precisely
\[
 (t,s,p)\in\{(5,7,5),(6,8,6),(5,7,4),(5,8,4),(5,9,3)\}.
\]

\emph{The five exceptions.}
First consider $(t,s,p)=(5,7,5),(6,8,6)$, so $p=t$ and
$h=t-1>\lfloor s/2\rfloor$. If a core edge is marked at one end of $P$, the arc restriction
forces every unused trace at the opposite end to be empty. At that opposite end, at most $h-1$ selected path
representatives and at most $\binom{p-1}k$ unused edges avoiding $C$
are available. Thus at least $D(n,r)-(h-1)-\binom{p-1}k=K-t+3\ge4$
core edges are marked at that end. Their tails cannot all
have distance at least $h$ from the first marked tail. Hence no
core edge is marked at either end, and each
end has at least four unused edges meeting $C$ in its assigned colour.

For $(t,s,p)=(5,7,4),(5,8,4)$, we have $h=3$, $r=4$ and $D(n,r)=11$.
A mark leaves at most two consecutive core vertices for an opposite
trace. Independence bounds the opposite unused edges by
$1+2\binom32=7$, counting $W$ and one of these core vertices.
With at most two path representatives, this forces an opposite mark.
The tail-distance restriction permits at most four opposite marks,
so each end has at least $11-2-4=5$ unused edges, at most one avoiding $C$.

In all four cases, each assigned-colour family of unused edges
meeting $C$ has at least four members. Form a bipartite graph with
these families as its vertex classes, joining distinct $r$-edges.
It is complete minus a matching, hence connected. Adjacent edges have distinct trace vertices at
distance at least $h+1$, so their traces share a singleton when
$2(h+1)>s$, or lie in one antipodal pair when equality holds.
Connectivity extends this to both families.

At $p=t$, the common vertex is a core neighbour with at least two
unused representatives. At $p=4$, a common singleton allows at most
$\binom43=4$ unused edges containing either endpoint and
contained in $W$ together with that singleton, contradicting the lower bound five.
The remaining possibility has $s=8$ and both antipodes in the trace
union at each end. The arc exclusion then permits no mark at either
end. At each end, the assigned colour has at least
$11-2=9$ unused edges containing that end.
Those avoiding one specified antipode lie in $W$ together with the
other antipode, giving at most four choices. Thus each antipode belongs to at least five unused edges
containing the relevant endpoint, giving two core neighbours with
at least two distinct representatives each.

Finally let $(t,s,p)=(5,9,3)$, so $h=2$ and $r=4$.
The first displayed sufficient condition holds, and the preceding
count gives
$11\le9+\binom{2+j}3-j.$
For $j=0,1,2$ the right sides are $9,9,11$, respectively.
Thus $j\le2$ forces $j=2$ and equality throughout: the trace union
has exactly $d+1$ incident core pairs, so is an interval. The unused
link at $x$ contains all four triples on $W\setminus\{x\}$ and
the two core neighbours with multiple representatives. Trace
independence makes these neighbours nonconsecutive.

Label the ordinary core cycle $C_9$ cyclically so these neighbours
are $0,i$, where $i=2,3,4$ after reversal if necessary. Each has at least two unused representatives; therefore any opposite trace
vertex distinct from it has distance at least three from it.
The possible opposite vertices are, respectively,
$\{5,6\}$, $\{0,3,6\}$ and $\{0,4\}$.
Colour $b$ has at least $11-1-c(z)\ge6$ unused edges at $z$,
since $c(z)\le\max\{0,6-d\}\le4$. A trace union of size at most
two permits only $\binom43=4$ such edges. Hence $i=3$, and all three vertices $0,3,6$
occur. Each belongs to at least two of those unused edges, since at
most four avoid any specified one of the three vertices. The reverse distance restriction confines the original trace union
to $\{0,3,6\}$, contradicting that it is an interval containing
$0,3$. Hence $j\ge3=t-p+1$.

If $p=1$, both core-neighbour sets have at least $t$ vertices.
Since $s<2t$, Lemma~\ref{lem:separation}(i), or the independence
number of the core cycle when the sets coincide, gives adjacent
vertices in the two sets. Their representatives are distinct by
trace independence, so they insert $x$ into $C$, a contradiction.
If $h\ge2$, distinct vertices in the two neighbour sets have core
distance at least $h+1$. The attachments can be chosen distinct: each has two choices for
$r\ge4$, while a triple cannot contain both ends of $P$ and two core vertices. If $p<q$, put $j=t-p+1\ge2$. Unequal
sets give $s\ge2j+2p=2t+2$ by Lemma~\ref{lem:separation}(ii). Equal
sets give $s\ge(p+1)j>2t+1-p\ge s$. Hence $p=q$ and $h=q-1$.
Write $A,B$ for these neighbour sets in the final detour.

\emph{Now suppose that $r=t$.} Here $q_0=q$, $W=V(P)$, and put $\kappa=r-h$.
If $h=1$, every unused edge in an unused colour has at least
$r-1$ core vertices. If some unused colour has no unused edge
containing $x$, then $d=c(x)\ge r+1$. By the two-mark exclusion, each unused edge contains at most
one head and at most one tail of a core edge marked at $x$. Thus the marked core edges form $c\ge1$ maximal proper paths and a
trace has size at most $s-d-c+2\le r-1$. Equality forces $s=2r-1$,
$q=3$ and $c=1$. The only possible trace is the vertex set of the unmarked path,
including its two ends. Since $q=3$, this gives at most three unused
edges in that colour, whereas degree summation gives at least
$\lceil nD(n,r)/r\rceil-s\ge q+3$.
Thus each unused colour has an unused edge containing $x$. Its
independent trace is incident with $2r-2$ distinct unmarked core
edges, so $c(x)\le1$. Choose distinct unused edges containing $x$ in two unused
colours; each colour offers at least $r$ choices. Their traces
are distinct $(r-1)$-sets and therefore have adjacent vertices
by Lemma~\ref{lem:separation}(i), yielding a longer Berge cycle.

Let $h\ge2$, and fix an unused colour. If $\kappa>0$ and this colour
has no unused edge at either end of $P$, then at least $\kappa+2$ core
edges are marked at each end. Choose that many tails
at each end. Distinct vertices from the two tail sets have core distance at
least $h$.
Equal sets give $s\ge h(\kappa+2)\ge2r$, and unequal sets give
$s\ge2(\kappa+2)+2(h-1)$ by \eqref{eq:cyclic-block-gaps}.
At $\kappa=0$ we have $q=h=r$ and $s=r+2$. Suppose neither end
has an unused edge of this colour meeting $V(C)$. Since $h$ exceeds
the diameter of the core cycle, the tail-distance restriction permits
only one core pair to be marked at either end, and that pair is marked
at both ends. Let $e$ be its representative, on the core pair
$a_0b_0$. The $r-1$ path representatives and $e$ are exactly the
$r$ selected representatives containing either endpoint $x,z$ of
$P$. Since $D(2r+2,r)=r+1$, each endpoint has at least one unused
incident hyperedge in the fixed unused colour. By the supposition,
it avoids $V(C)$, and hence it must be $W=V(P)$, the only $r$-set
contained in $V(P)$. Thus $W$ is an unused hyperedge. Replace
one path representative by $W$ in the unused colour, preserving the
core of $P$, and let $g$ be the released representative. Then
$g\ne W$, contains both $x,z$, and meets $V(C)$; choose
$d\in g\cap V(C)$. There are two arcs from $d$ to the marked pair
$a_0b_0$ that traverse that pair, one ending at $a_0$ and the other
at $b_0$. Their lengths sum to $s+1$. Choose the shorter one, let
$a$ be its terminal endpoint and $b$ the other endpoint of the marked
pair. Its length satisfies
$\rho\le\lfloor(s+1)/2\rfloor=\lfloor(h+3)/2\rfloor\le h$.
Replace this arc by the detour $d,z,P,x,a$, using $g$ on $dz$, the
representatives of $P$ from $z$ to $x$, and $e$ on $xa$. The detour
has $h+1$ edges, so it is strictly longer than the removed arc. The new cycle is rainbow and uses distinct representatives:
$g$ has the released path colour, $e$ has a cycle colour, and $W$
uses the previously unused colour; all other representatives retain
their original distinct assignments. This contradicts the maximality
of $C$.

A nonempty trace at one end permits a core edge marked at the
opposite end only in a cyclic gap of length at least $2h+1$. If
$\kappa=0$, then $s=h+2<2h+1$, so no such gap exists. Otherwise,
choose an unused edge at the first end with a nonempty trace.
This trace contains at least $r-h=\kappa$ core vertices and is
independent, so it contains an independent $\kappa$-set. Apply Lemma~\ref{lem:cyclic-gap}; the
largest gap has length at most
$s-2(\kappa-1)=2h+3+\varepsilon-q$. Thus a marked edge forces
$q=3$, $\varepsilon=1$, and equality in the gap bound, with every
other gap of length two. Only the central core edge of the long gap
can be marked.
Each unused colour therefore has at least $\kappa+1$ unused edges
at the opposite end, and at least two when $\kappa=0$. At least one
of these edges meets $V(C)$, so the same reasoning applies in
the other direction. If $h<q$, then $\kappa>0$ and two colours remain. Distinct unused edges
at opposite ends in different colours have distinct $\kappa$-subsets in
their traces: only one edge fits in $W$ and $\kappa$ core vertices.
Equation~\eqref{eq:cyclic-block-gaps} with $\ell=h$ gives
$s\ge2\kappa+2h=2r$. Thus $h=q$, with one unused colour $a$.

Orient $P=x_0\ldots x_{q-1}$ with $x=x_0$. If a terminal representative does not contain
$W$, orient $P$ with that representative last, and let
$B$ be its trace. Otherwise both terminal representatives contain
$W$ and may be interchanged with their assigned
colours; let $B$ be the union of their traces.
Let $A$ be the union of traces of unused colour-$a$ edges
containing $x$. Both sets have at least $\kappa+1$ vertices when $\kappa>0$: a terminal
avoiding part of $W$ has that many core vertices, and distinct edges
cannot both lie in $W$ and $\kappa$ core vertices. At $\kappa=0$ both are
nonempty, since at most one terminal equals $W$.
Attach the subpath with core order $x_0,\ldots,x_{q-2}$ to the
cycle using an unused colour-$a$ edge at the first end and the
reserved last representative at the other. This shows that
distinct vertices in $A$ and $B$ have core distance at least $q$.

\emph{The final detour.} In both cases, $W=V\setminus V(C)$.
Use the first $m=q-1$ vertices $x_0=x,\ldots,x_{m-1}$ of $P$.
For $c\in B$, reserve an attachment on $x_{m-1}c$: unused in
colour $b$ when $r<t$, or terminal when $r=t$.
Write $j=t-q$. We have $|A|,|B|\ge j+1$ and distinct cross vertices
have distance at least $q$. At $j=0$, $s=q+2<2q$ forces a common
singleton. For $j>0$, unequal sets give $s\ge2(j+1)+2(q-1)=2t$
by Lemma~\ref{lem:separation}(ii). Hence $A=B=S$, with gaps of
length at least $q$. Since $(j+2)q>s=q+2j+1+\varepsilon$, $|S|=j+1$.
The largest gap is $q+2$ at $j=0$ and at most
$s-jq=q+1+\varepsilon-(q-2)j\le4$ otherwise.

A core edge marked at $x$ requires a gap in $S$ of length at
least $2m+1$.
From either end of the marked pair, none of the first $m$ vertices
reached across that pair can belong to $S$. Otherwise the mark, the
path on $x_0,\ldots,x_{m-1}$ and the last attachment would
replace an arc of at most $m$ edges by $m+1$. The gap bounds exclude
this except at $q=t=3$. In that case $\kappa=0$, and a nonempty opposite trace, together with
the entire three-vertex path $P$, excludes a marked core edge.
Thus $c(x)=0$. Every core vertex is within distance $m$ of $S$,
also when $q=3$ and the ordinary core is $C_5$.
Put $E=W\cup S$, so $|E|=t+1$.
Every unused colour-$a$ edge at $x$ lies in $E$, by the
definition of $A$ at $r=t$. At $r<t$ an attachment outside $S$ would
give a longer detour using a distinct opposite representative.

Degree supplies $g\in H_a(x)$ leaving $E$. It is selected and,
since $c(x)=0$, represents a pair $x_ix_{i+1}$ of $P$.
Choose $d\in g\setminus E$, a nearest vertex $c\in S$, and a representative
on $x_{m-1}c$. At $r=t$, every terminal representative that can be placed last lies
in $E$, so $g$ keeps its pair and $i\le m-2$.
At most $h-1$ selected edges contain $x$. Thus at least
\[
 D(n,r)-(h-1)-\binom{t-1}{r-1}=\binom{t-1}{r-2}+2-h
\]
unused colour-$a$ edges contain $x,x_{i+1}$. For $r<t$ there are
at least two; for $r=t$ the fixed last attachment is selected.
In either case, choose a representative distinct from that attachment.

Replace a shortest $d,c$-arc of the core cycle by
\[
 d,x_i,x_{i-1},\ldots,x_0,x_{i+1},\ldots,x_{m-1},c,
\]
omitting empty segments. Move $g$ from its old pair to $dx_i$ in its original colour; use
the chosen colour-$a$ edge on $x_0x_{i+1}$ and the fixed attachment
on $x_{m-1}c$. All other path representatives and colours remain. The attachments
are unused or reserved, $g$ occurs once, and the inserted vertices
are distinct and outside $C$. This gives a rainbow Berge cycle of
order $s-\rho+m+1>s$, where $\rho\le m$ is the removed arc length.
\end{proof}

\section{Berge Hamilton cycles from shadow graphs}\label{sec:binomial}

Let $\mathcal F=(F_1,\ldots,F_n)$ be the shadow graph system of
$\mathcal H$, as defined in Section~\ref{sec:notation}.

\begin{lemma}\label{merge:lifting}
Let $\mathcal H=(H_1,\ldots,H_n)$ be an $r$-graph system on $n$ vertices,
$3\le r<n$, with shadow graph system $\mathcal F$. Let $C$ be a
rainbow ordinary cycle of length greater than $r$ in $\mathcal F$.
Prescribe distinct representatives for $d$ of its pairs. These
extend to a lifting if every remaining pair has codegree at least
\[
 \begin{cases}
 r-1,&d\le1,\\
 r,&d=2,\\
 r+d-3,&d\ge3
 \end{cases}
\]
in its assigned layer. Every rainbow ordinary forest whose pairs
have codegree at least $r-1$ in their assigned layers admits distinct
representatives.
\end{lemma}
\begin{proof}
Let $S(e)$ be the allowed representative set of each coloured pair
$e$, restricted to the chosen singleton for a prescribed pair.
Hall's theorem gives distinct representatives exactly when
$|\bigcup_{e\in\mathcal A}S(e)|\ge|\mathcal A|$ for every subfamily
$\mathcal A$. For a forest, each pair has at least $r-1$ representatives,
and an $r$-edge contains at most $r-1$ of its pairs. Hence every subfamily $\mathcal A$ satisfies
$(r-1)|\mathcal A|\le(r-1)|\bigcup_{e\in\mathcal A}S(e)|$.

A Hall-deficient subfamily of cycle pairs has $a\ge1$ unprescribed
pairs, $b\le d$ prescribed pairs, and $\tau\le u\le a+b-1$
representatives in their union, where $\tau$ is the stated bound. As $|V(C)|>r$, each representative
contains at most $r-1$ cycle pairs, so
\begin{equation}\label{eq:lifting-strong-incidence}
 \tau a\le(r-1)u.
\end{equation}
Including the prescribed representatives gives
$\tau a+b\le(r-1)(a+b-1)$.
For $d\le1$, this requires $(r-2)b\ge r-1$, impossible.
For $d=2$, it gives $a\le(r-2)b-(r-1)\le r-3$, whereas
$r\le u\le a+b-1\le r-2$, a contradiction.

Let $d\ge3$, so $\tau=r+d-3$. Equation~\ref{eq:lifting-strong-incidence}
excludes $b=0$, and $u\ge\tau$ gives $a\ge\tau-b+1\ge r-2$.
Write $j_1,\ldots,j_q$ for the lengths of the maximal subpaths of
selected unprescribed pairs. Only the vertex sets of consecutive $(r-1)$-edge subpaths contain
$r-1$ selected unprescribed pairs. Thus at most
$\sum_{s=1}^q\max\{0,j_s-r+2\}\le a-r+2$
representatives contain $r-1$ selected unprescribed pairs; the others
contain at most $r-2$. Consequently,
\[
 \tau a\le(r-2)u+a-r+2\le(r-1)a+(r-2)(b-2),
\]
so $(d-2)a\le(r-2)(b-2)\le(r-2)(d-2)$.
It follows that $a=r-2$, $b=d$, $u=\tau$, with equality throughout.

Every representative contains all $r-2$ unprescribed pairs.
These pairs span at most $r$ vertices; equality would permit only
one representative, contrary to $u=\tau\ge r$. Thus they form an $(r-1)$-vertex path.
A representative adds one vertex to this path; its only additional
cycle pairs are at the two ends, contradicting $b=d\ge3$.
\end{proof}

For the remainder of this section, recall that $\varepsilon\in\{0,1\}$, put $k=r-1\ge3$, and assume
\begin{equation}\label{eq:shadow-range}
 2k\le t+\varepsilon,\qquad
 t\ge
 \begin{cases}
 13,& k=3,\\
 9,& k\ge4.
 \end{cases}
\end{equation}
Define the \emph{strong shadows} by
$E(G_i)=\{uv:|S_i(uv)|\ge r+\varepsilon\}$.
A colour-$i$ pair is \emph{strong} if it belongs to $G_i$, and
\emph{weak} if it belongs to $F_i$ but not $G_i$. An ordinary
rainbow path, cycle or forest is strong if all its pairs are strong.

\begin{lemma}\label{lem:shadow-degree}
Let $\mathcal H=(H_1,\ldots,H_n)$ be an $r$-graph system on
$n=2t+1+\varepsilon$ vertices, $\varepsilon\in\{0,1\}$.
Put $k=r-1\ge3$ and suppose $2k\le t+\varepsilon$, with
$t\ge13$ if $k=3$ and $t\ge9$ otherwise. Assume
$\delta_1(H_i)\ge D(n,r)$ for every $i$.
Define $E(G_i)=\{uv:|S_i(uv)|\ge r+\varepsilon\}$.
Then $\delta(F_i)\ge t+1$ and $\delta(G_i)\ge t$ for every $i$.
If $n$ is even and $\mathcal H$ has no rainbow Berge Hamilton cycle,
the $G_i$ do not contain a common $K_{t,t+2}$.
\end{lemma}
\begin{proof}
Fix $i,v$. If $d_{F_i}(v)\le t$, its link has at most $\binom tk$
edges, contrary to degree. Thus $\delta(F_i)\ge t+1$.

Suppose $d_{G_i}(v)\le t-1$. In the link of $v$, at most
$t-1$ vertices have degree greater than $k+\varepsilon$.
Lemma~\ref{lem:link-incidence}, with $N=2t+\varepsilon$,
$s_0=t-1$ and $\beta=k+\varepsilon$, gives
$d_{H_i}(v)\le\binom{t-1}k+(k+\varepsilon)(t+1+\varepsilon).$
It suffices that $\binom{t-1}{k-1}\ge(k+\varepsilon)(t+1+\varepsilon)$.
For $k=3$, use $\binom{t-1}2\ge4(t+2)$ for $t\ge13$.
For $4\le k\le(t+1)/2$, the ratio of consecutive terms of
$\binom{t-1}{k-1}/(k+1)$ is $(t-k)(k+1)/(k(k+2))\ge1$
until the last term. Hence, for $t\ge9$,
\[
 \frac{\binom{t-1}{k-1}}{k+1}\ge\frac{\binom{t-1}3}{5}\ge t+2,
\]
which proves $\delta(G_i)\ge t$.

For the last assertion, suppose every pair joining $X$ to $Y$
is strong in every colour, where $|X|=t$ and $|Y|=t+2$.
Since $\delta(F_i)\ge t+1$, each $F_i[Y]$ has minimum degree at
least one. Choose an edge of $F_i[Y]$ and a representative $f\in H_i$.
In another colour $j$, choose a pair of $F_j[Y]$ meeting
$Y\setminus f$; both the pair and its representative differ from
the first choices. These pairs and
the isolated vertices of $Y$ form $t$ path components. Alternate
them with $X$ around a Hamilton cycle, assigning the remaining colours
to the strong cross pairs. The two possible weak pairs have distinct
representatives, so Lemma~\ref{merge:lifting} gives a contradiction.
\end{proof}

\subsection{An ordinary shadow cycle omitting one vertex}\label{sec:odd}
\begin{proposition}\label{prop:even-one}\label{sec:even-one}
Let $\mathcal H=(H_1,\ldots,H_n)$ be an $r$-graph system on
$n=2t+1+\varepsilon$ vertices, $\varepsilon\in\{0,1\}$.
Put $k=r-1\ge3$ and suppose $2k\le t+\varepsilon$, with
$t\ge13$ if $k=3$ and $t\ge9$ otherwise. Assume
$\delta_1(H_i)\ge D(n,r)$ for every $i$.
If the shadow graph system of $\mathcal H$ has a rainbow ordinary
$(n-1)$-cycle with at most $\varepsilon$ pairs of codegree below
$r+\varepsilon$, each of codegree at least two in its assigned layer,
then $\mathcal H$ has a rainbow Berge Hamilton cycle.
\end{proposition}
\begin{proof}
Suppose no rainbow Berge Hamilton cycle exists. Write the ordinary
cycle in core order $(x_1,\ldots,x_m)$, with omitted vertex $y$
and $m=2t+\varepsilon$. Relabel colours so that $0$ is
unused and $e_h=x_hx_{h+1}:h$; subscripts are modulo $m$.
If a weak pair $e_\ell$ exists, reserve two representatives
$\mathcal U\subseteq S_\ell(x_\ell x_{\ell+1})$; otherwise put
$\mathcal U=\varnothing$. Every Hamilton cycle constructed below has at most
$3+\varepsilon$ weak pairs. Lemma~\ref{merge:lifting} reduces the
proof to assigning distinct representatives to the new weak pairs
while retaining a reserve for $e_\ell$ if it survives.

For $i\in[m]$, put
\[
 L_i=\{h\in[m]:x_hx_i\in E(F_h),\ h+1\ne i\},
 \qquad d_i=|L_i|.
\]
Also put
\[
 J=\{j\in[m]:yx_{j+1}\in E(G_0)\},
 \qquad \sigma=|\{h\in[m]:yx_h\in E(F_h)\}|.
\]
Then
\begin{equation}\label{eq:rotation-demand}
 \sigma+\sum_i d_i=mt+\rho,\qquad
 \rho=\sum_h\bigl(d_{F_h}(x_h)-t-1\bigr)\ge0 ,
\end{equation}
by Lemma~\ref{lem:shadow-degree} and counting the neighbours of
each $x_h$ other than $x_{h+1}$.
Delete $e_i$ and attach $y$ at $x_{i+1}$, or insert it at $e_j$,
$j\ne i$. Orient the resulting ordinary path towards $x_i$,
preserving all retained edge orientations. Replacing a retained
$e_h$, $h\in L_i$, by $x_ix_h:h$ makes $x_{h+1}$ an endpoint.
Including $x_i$, this gives $1+d_i$ endpoints after attachment and
$d_i+\mathbf1_{j\notin L_i}$ after insertion.

\emph{Reserve conflicts.} Two exchanges suffice. First, suppose insertion
at $e_h$ uses $f_0:0$ on $yx_{h+1}$ and $f_h:h$ on $yx_h$,
exhausting the reserves while $e_\ell$ survives. Both edges contain
$y,x_\ell,x_{\ell+1}$ and belong to $H_\ell$. Instead insert $y$
at $e_\ell$, using $f_0:0$ on $yx_{\ell+1}$ and $f_h:\ell$ on
$yx_\ell$, thereby removing the reserved pair.

Second, let $z$ be the fixed end after attachment or insertion.
Suppose rotation on $x_ix_h$ with $g_1:h$ and closure on $zx_{h+1}$
with $g_2:j$ exhaust the reserves while $e_\ell$ survives, and
$zx_i\notin E(F_j)$. Then $x_i\notin g_2$, so $x_i$ is not an end
of $e_\ell$. Rotate at $e_\ell$ instead, using $g_1:\ell$ on
$x_ix_\ell$ and $g_2:j$ on $zx_{\ell+1}$. The path orientation gives a cycle on the same vertices, retaining
the attachments at $y$, restoring the strong pair $e_h$ in colour
$h$, and removing $e_\ell$.
The other weak representatives remain distinct.

With strong attachments at $y$, no attainable endpoint permits
closure in the unused shadow graph. Direct closure uses at most
one reserve. A representative shared by closure and rotation would
also close the unrotated path on $zx_i$. Distinct representatives exhausting the
reserves are handled by the second exchange. Every closure would
therefore lift.

\emph{Endpoint bounds.} Put $B=J+1$, the cyclic shift of $J$,
and $K=N_{F_0}(y)$.
If $yx_h\in E(F_h)$ and $yx_{h+1}\in E(F_0)$, insertion at
$e_h$, together with the first exchange above, forces the two
pairs to have the same unique representative. Interchanging its assigned colour shows that neither pair is strong
in colour $0$. There are at most $m-|K|$
indices with $x_{h+1}\notin K$. Each remaining index counted by $\sigma$ gives a core pair in
$K\setminus\{x_i:i\in B\}$. This proper subset induces at most
$\max(0,|K|-|J|-1)$ edges of $C$, so
$\sigma\le m-|K|+\max(0,|K|-|J|-1)\le t+\varepsilon-1$.
For $j\in J$, delete $e_j$ and attach $yx_{j+1}:0$. The
$1+d_j$ attainable endpoints avoid $N_{F_j}(y)$. Since
$|N_{F_j}(y)|\ge t+1$, we obtain $d_j\le t+\varepsilon-2$.
In even order, equality makes $N_{F_j}(y)$ the complementary set,
so $yx_{j+1}\in E(F_j)$. Recolour the attachment with $j$.
A representative on $yx_j:0$ differs from it because
$yx_j\notin E(F_j)$; the first exchange gives
\begin{equation}\label{eq:prec-exclusion}
 \varepsilon=1,\ j\in J,\ d_j=t-1\quad\Longrightarrow\quad yx_j\notin E(F_0) .
\end{equation}

\smallskip\emph{The degree comparison.} If $d_i\le t$ for every $i$,
then \eqref{eq:rotation-demand}
and the bound on $d_j$ for $j\in J$ give
$0\le\rho\le\sigma-(2-\varepsilon)|J|.$
Together with $|J|\ge t$ and $\sigma\le t+\varepsilon-1$,
this forces $\varepsilon=1$, $|J|=t$ and $d_j=t-1$ for every
$j\in J$. We prove the same conclusion if some $d_i>t$.

Now fix $i$ with $d_i>t$ and put $v=x_{i+1}$. Consider
$e_j=aw$, in either orientation, with
$ya\in E(F_i)$ and $yw\in E(G_0)$. Insertion at $e_j$ excludes
$ya\in E(F_j)$, and therefore $i\ne j$. Cut $e_i$ and insert $y$ at $e_j$
using $ya:i,yw:0$. Among the possible endpoints
$T_{i,j}=\{x_i\}\cup\{x_{h+1}:h\in L_i\setminus\{j\}\}$,
at least $d_i-t$ are adjacent to $v$ in $F_j$, since $v$ has at
most $t$ nonneighbours other than itself. The exclusion for
strong attachments gives $ya\notin E(G_i)$. The indices of $N_{G_i}(y)$ avoid
$J\cup(J+2)$. Both sets $J$ and $N_{G_i}(y)$ have size at least $t$.
For $m=2t$ disjointness gives $|J|=t$ and $J+2=J$, so $J$ is a parity class. For $m=2t+1$, shifting by two is transitive, so $|J|=t$ and
$|(J+2)\setminus J|=1$.
In both cases label $J=\{1,3,\ldots,2t-1\}$,
$B=\{2,4,\ldots,2t\}$, with $N_{G_i}(y)=\{x_h:h\in B\}$.
When $n$ is even, $m$ is the only index outside $J\cup B$.

Choose $a=x_j\in N_{F_i}(y)$ with $j\in J$. When $n$ is even
and $yx_m\in E(F_i)$, instead choose $a=x_m$ and $j=m-1$.
In either case, write $e_j=aw$ with $yw\in E(G_0)$.
If $ya:i$ uses a reserve while $e_\ell$ survives, neither end of
$e_\ell$ lies in $B$, or a strong colour-$0$ pair inserts $y$.
Thus $\varepsilon=1$, $\ell=m$, $a=x_m$ and $i\in B$. Whenever $i\in B$ and $a=x_m$,
we have $vx_m\notin E(F_j)$, or the arcs of $C-e_i-e_j$ close via
$x_i y:i,yx_{m-1}:0,vx_m:j$ with only the last new pair weak.
Every closing representative then avoids $x_m$, hence both reserves. Direct closure leaves a
reserve because its representative differs from the attachment ($ya\notin E(F_j)$).
Thus $vx_i\notin E(F_j)$.

For $h\in L_i\setminus\{j\}$ with $vx_{h+1}\in E(F_j)$,
choose representatives $f\in H_i$, $g_1\in H_h$ and $g_2\in H_j$
for the pairs $ya$, $x_ix_h$ and $vx_{h+1}$, respectively.
Replace $e_i,e_j,e_h$ by these pairs and the strong pair $yw:0$
to close. The representative $g_2$ differs from $f$ and $g_1$: equality would
put $ya$ or $vx_i$ in $F_j$.

If $f\ne g_1$, all three are distinct. A reserve conflict cannot
involve $f,g_2$, since a reserved attachment makes $g_2$ avoid both
reserves; a conflict involving $g_1,g_2$ is resolved by the second
exchange. Thus only $f,g_1$ as the two reserves remain, forcing
$i\in B$, $a=x_\ell=x_m$, $\ell=m$.
If $h\in J$, delete $e_m,e_h$ from the original cycle and insert
\[
 x_1y:m\text{ using }f,\qquad
 yx_{h+1}:0,\qquad x_hx_m:h\text{ using }g_1 .
\]
The retained arcs form one cycle. The middle pair is strong, the
other two have distinct representatives, and $e_m$ is removed.
Thus the cycle lifts, so $h\notin J$.

If no distinct $f,g_1$ can be chosen, both candidate sets are
$\{g_1\}$, with $y,a,x_i,x_h\in g_1$. If $h\in J$, insert $y$ at $e_h$ using
$g_1$ on $yx_h:h$ and the strong pair $yx_{h+1}:0$. Lemma~\ref{merge:lifting} supplies the strong representative while
leaving a reserve when needed. Again $h\notin J$.

In both obstructions $a\in g_1$, and $g_1,g_2$ leave a reserve: discard
the attachment reserve $f$ in the first; in the singleton case,
$g_1\in\mathcal U$ forces $g_2\notin\mathcal U$.
If $i,h\in B$, then $a\ne x_i,x_h$. Delete the disjoint pairs
$ya,x_ix_h$ from the proposed cycle.
Exactly one of $\{yx_i,ax_h\}$ and $\{yx_h,ax_i\}$ reconnects
the resulting paths. Use colour $i$ on the pair at $y$, which is
strong because $i,h\in B$, and $g_1:h$ on the pair at $a$.
Discard $f$, retaining $g_2$ and $yw:0$. The weak representatives
remain distinct and leave a reserve, so the cycle lifts.
Hence $i,h$ cannot both lie in $B$.
Since $i,h\notin J$, no rotation permits closure when $n$ is odd.

When $n$ is even and $i\in B$, the only rotation index
permitting closure is $m$.
Then $a=x_m$ by our choice, since the common representative $g_1$ also gives
$yx_m\in E(F_i)$. Thus $j=m-1$. At least $t$ vertices in $T_{i,j}\setminus\{x_1\}$ cannot close.
Together with the additional nonneighbour $x_m\ne v$, they
contradict $d_{F_j}(v)\ge t+1$.

It remains to consider $i=m$, with $d_m=t+\eta>t$. In this case,
the obstruction using both reserved representatives is impossible.
There are at least $\eta$ distinct indices $h\in B$ at which a
rotation permits closure, and each satisfies
$yx_h\in E(F_h)$ and $|S_h(x_mx_h)|=1$.
Also $h\ne m-1$ by the definition of $L_m$. If $d_{F_h}(x_h)=t+1$, at most one colour-$h$ edge at $x_h$
contains $x_m$. Degree therefore forces all $\binom tk$ link sets
on its other $t$ neighbours. These include $y$ (via $g_1$) and
$x_{h+1}$ (via the original core pair), so one such link set
contains both.
Use such a representative on $yx_{h+1}:h$ and the strong pair
$yx_h:0$ to insert $y$ at $e_h$. Thus $\rho\ge\eta$, while
$\sigma+\sum_i d_i\le t+t(t-1)+t^2+t+\eta=mt+\eta$.
Equality in \eqref{eq:rotation-demand} gives $d_j=t-1$ on $J$.

Thus $n$ is even, $|J|=t$ and $d_j=t-1$
for every $j\in J$. Apply \eqref{eq:prec-exclusion} in both
orientations of $C$. Every core neighbour of $S=N_{G_0}(y)$ lies
outside $K$. However, $|V(C)\setminus K|\le t$, whereas
$|N_C(S)|=t+|(S+2)\setminus S|>t.$
Here $S+2$ denotes the cyclic shift by two, and the strict inequality
holds because $|V(C)|=2t+1$ and $|S|=t$. This contradiction proves
the proposition.
\end{proof}

\subsection{Long cycles in the shadow graph system}\label{sec:even-long}
\begin{proposition}\label{prop:shadow-conclusion}
Let $\mathcal H=(H_1,\ldots,H_n)$ be an $r$-graph system on
$n=2t+1+\varepsilon$ vertices, $\varepsilon\in\{0,1\}$.
Put $k=r-1\ge3$ and suppose $2k\le t+\varepsilon$, with
$t\ge13$ if $k=3$ and $t\ge9$ otherwise. Assume
$\delta_1(H_i)\ge D(n,r)$ for every $i$.
Then $\mathcal H$ has a rainbow Berge Hamilton cycle.
\end{proposition}
\begin{proof}
Suppose no rainbow Berge Hamilton cycle exists.
Lemma~\ref{lem:shadow-degree} gives $\delta(F_i)\ge t+1$ and
$\delta(G_i)\ge t$.

\emph{Odd order.} If $n=2t+1$, \cite[Lemma~3.3]{LLL} gives a strong
rainbow ordinary cycle of length at least $n-1$, unless every
strong shadow consists of the same two $(t+1)$-cliques sharing a vertex $z$. In the exception,
$\delta(F_i)\ge t+1$ supplies a pair joining the cliques away from
$z$. Join its ends to $z$ by spanning paths in the two cliques and
assign the remaining colours to their strong pairs. The resulting
Hamilton cycle lifts by Lemma~\ref{merge:lifting}. Otherwise, lift
a strong Hamilton cycle directly or apply
Proposition~\ref{prop:even-one} to a strong $(n-1)$-cycle.

\smallskip\noindent\emph{Even order.}
Let $n=2t+2$. By \cite[Theorem~1.3]{LLL}, there is a strong
rainbow ordinary spanning path unless all strong shadows are the
same union of two disjoint cliques or contain a common $K_{t,t+2}$.
Lemma~\ref{lem:shadow-degree} excludes the latter.

In the clique exception, both cliques $A,B$ have order $t+1$.
Choose a cross pair in $F_i$ with representative $f$, relabelling
so $|f\cap B|\ge2$. Choose $u\in A\setminus f$, possible since
$r<t+1$, and a cross neighbour $v\in N_{F_j}(u)\cap B$ in
colour $j\ne i$. For $p\in A\cap f$ and
$q\in(B\cap f)\setminus\{v\}$, the disjoint pairs $pq:i,uv:j$
have distinct representatives because $u\notin f$. Join their ends
by spanning paths in $A,B$, assigning the remaining colours to the
strong pairs. This Hamilton cycle lifts by Lemma~\ref{merge:lifting}.

We first find a strong cycle of length at least $t+2$.
Split a strong rainbow ordinary spanning path into consecutive
blocks $A,B$ of order $t+1$, with bridge $a^*b^*:\kappa$ and unused
colour $c$, and suppose no such cycle exists.
Write the spanning path in the order
$(a_t,\ldots,a_1,a_0=a^*,b_0=b^*,b_1,\ldots,b_t)$.
Its endpoints are $a_t,b_t$. A neighbour of either endpoint outside
its own block would close a strong rainbow cycle of length at least
$t+2$. The degree bound therefore gives
$N_{G_c}(b_t)=\{b_0,\ldots,b_{t-1}\}$ and
$N_{G_c}(a_t)=\{a_0,\ldots,a_{t-1}\}$. The same conclusion holds
after each rotation retaining the bridge, for the new endpoint and
unused colour. Write $b_ib_{i+1}:\gamma_i$ for the successive
coloured pairs of the block path.
For $1\le i<t$, rotate first on $b_ib_t:c$, then on
$b_{i-1}b_{i+1}:\gamma_i$. The successive endpoint and unused-colour
pairs are $(b_{i+1},\gamma_i)$ and $(b_i,\gamma_{i-1})$.
Recolouring $b_ib_t$ with $\gamma_{i-1}$ frees $c$ again; for
$i=t-1$ the first step only recolours a pair. Apply the same rotations independently in $A$. This prescribes
any ends $a\in A\setminus\{a^*\}$, $b\in B\setminus\{b^*\}$
while retaining the bridge and unused colour $c$.

Closure and lifting exclude $ab\in E(F_c)$, so
$\delta(F_c)\ge t+1$ forces $ab^*,a^*b\in E(F_c)$.
A bridge representative contains another vertex; interchange blocks
so it can serve as $b$. Replacing the bridge by
$a^*b:\kappa,ab^*:c$ forms a Hamilton cycle with distinct new
representatives, or $ab\in E(F_c)$. Its at most two weak pairs
allow lifting. Hence a strong cycle of length at least $t+2$ exists.

Next let $C$ be a longest strong rainbow ordinary cycle, of
length $s$. We prove $s\ge2t$. Suppose instead that $q=n-s\ge3$. Let $P$ be a longest strong rainbow ordinary path outside $C$ in colours
unused by $C$, of order $p\le q\le t$. Delete an endpoint if
$p=q$, obtaining a path of order $\ell=p-\mathbf1_{p=q}$ with
two unused colours. Let $S_a,S_b$ be the sets of strong neighbours
on $C$ of its two ends in two distinct such colours. By maximality
if $p<q$, and by the degree condition otherwise,
$|S_a|,|S_b|\ge t-p+1$. If $x\in S_a$ and $y\in S_b$ are
 distinct, then $\operatorname{dist}_C(x,y)\ge\ell+1$; a shorter
arc would be replaced by the corresponding detour through $P$,
lengthening $C$. We have $s\ge t+2$ from the preceding paragraph
and $s=n-q\le2t+2-q$. Lemma~\ref{lem:separation}(iv), applied
to $S_a,S_b$, gives $p=q=t$ and $s=t+2$.
Return to the original $p$-vertex path $P$ before the endpoint
was deleted. Put $X=V(P)$, let $T$ be the colours unused by $C$, and set
$h=\lfloor t/2\rfloor$. A subpath of length $h$ or $h+1$ leaves three $T$-colours unused.
Its ends have core neighbours in each; opposite neighbours in
different colours coincide, or the subpath replaces a shorter arc.
The third colour makes all these neighbourhoods the same singleton.

The graph on the indices $[t]$ of $P$ with steps $h,h+1$ is
connected: $i,i+1$ share $i+h+1$ for $i<t-h$, and larger indices
move back by $h$. Thus in the colour $c\in T$ unused by $P$, all
vertices of $X$ have the same unique core neighbour. Degree $t$
makes $G_c[X]$ complete. Starting with $c$ unused, recolour one edge
of $P$ with the currently unused colour. The former colour of that edge
is now unused; use it to recolour the next edge. Continue one edge at
a time. At every stage $P$ remains rainbow, so the same argument applies
to the newly unused colour. After all $t-1$ path edges have been treated,
the property holds for every $T$-colour. A path of length $h$ in $X$
avoiding two specified colours shows that their singleton is the same
vertex $x$.

Put $Y=V(C)\setminus\{x\}$, and let
$xz_1:\alpha_1$ and $xz_2:\alpha_2$ be the two pairs of $C$
incident with $x$. Suppose a colour $i$ used by $C$ has a
strong pair $wz:i$ with $w\in X$ and $z$ an internal vertex of
$C-x$. Follow the $xz$-arc avoiding colour $i$, of length at least two,
then $zw:i$, all of $X$ in $T$-colours, and return to $x$.
This longer strong cycle is forbidden. Thus each $C$-colour gives
minimum degree at least $t-3$ in $X$, since a vertex has at most
three strong neighbours outside. Each $T$-colour gives minimum
degree at least $t-1$ in $Y$, since no strong pair joins $X,Y$.

Choose the colour $d$ of an edge of $C-x$ incident with neither
$z_1$ nor $z_2$. Since $\delta(F_d)\ge t+1$, there is a shadow
pair $wz:d$ with $w\in X$ and $z\in Y$. Choose an edge $zz':e$
of $C-x$ with $e\ne d$, and choose $j\in\{1,2\}$ such that
$z_j\notin\{z,z'\}$. By \cite[Theorem~4]{SWW}, span $Y\setminus\{z\}$ from $z'$ to
$z_j$ in $t-1$ chosen $T$-colours, and span $X$ from $w$ to any
$w'\ne w$ in the colours of $C$ outside $\{d,e,\alpha_j\}$.
Each system has $t$ vertices and $t-1$ colours; its minimum degrees
are at least $t-2$ and $t-3$, respectively, both at least $(t+1)/2$
for $t\ge9$. Join the paths with $z'z:e,zw:d,xz_j:\alpha_j$ and
$w'x$ in the one remaining $T$-colour. All colours are distinct and only $zw$ can be weak, so
Lemma~\ref{merge:lifting} gives a contradiction. Thus $s\ge2t$.

Finally, a strong ordinary Hamilton cycle lifts, and
Proposition~\ref{prop:even-one} excludes a strong ordinary
$(n-1)$-cycle. Hence $C$ has length $2t$; let $u,v$ be its omitted
vertices and $c,d$ its unused colours.
For $w\in\{u,v\}$ and $i\in\{c,d\}$, define
\[
 B_i(w)=N_{G_i}(w)\cap V(C),\qquad
 A_i(w)=\{z\in V(C):|S_i(wz)|\ge2\}.
\]
Then $B_i(w)\subseteq A_i(w)$ and $|B_i(w)|\ge t-1$.
No core edge joins $A_c(w)$ to $B_d(w)$, or insertion of $w$
gives an ordinary $(n-1)$-cycle covered by
Proposition~\ref{prop:even-one}. Interchange $c,d$ for the reverse exclusion.

Index $C$ cyclically modulo $2t$. For $X\subseteq V(C)$,
let $X+2$ be its shift by two.
Then $|N_C(X)|=|X|+|(X+2)\setminus X|$.
Each partly occupied parity class contributes to the difference,
so $|X|\ge t-1$ implies $|N_C(X)|\ge t$, with equality only
for a parity class missing at most one vertex.
The adjacency exclusion gives $|A_i(w)|\le t$; equality makes
$A_i(w)$ a parity class and confines the other $A$-set to it.
In particular, $uv\notin E(G_i)$ forces
$A_i(w)=B_i(w)$ of size $t$ for both $w=u,v$.

Suppose both $A$-sets at $u$ are contained in a parity class
$X$, and both $A$-sets at $v$ are contained in the opposite class
$Y$. At most $t$ edges of $H_c(u)$ meet $Y$, since each pair $uz$
with $z\in Y$ has codegree at most one. Among the remaining edges,
those containing $v$ give pairwise disjoint $(k-1)$-subsets of $X$:
otherwise some pair $vz$, $z\in X$, would have codegree at least
two. There are at most $t/(k-1)$ such edges. Thus fewer than $t+t/(k-1)\le3t/2$ of the $k$-subsets of $X$
are absent from the link of $u$ in $H_c$.

A nonedge of $G_c[X]$ belongs to at least
\[
 \binom{t-2}{k-2}-(k+1)\ge\frac43\binom k2
\]
of these missing $k$-sets. Divide by
$\binom k2$. The binomial term becomes
$\binom tk/\binom t2\ge(t-2)/3$, while
$(k+1)/\binom k2\le4/3$. For $t=9$, where $k\ge4$, the latter
bound improves to $5/6$; these estimates cover
\eqref{eq:shadow-range}. Counting incidences between nonedges and
missing $k$-sets shows that $G_c[X]$ has fewer than $9t/8$
nonedges. Interchanging $u,c,X$ with $v,d,Y$ gives the same bound
for $G_d[Y]$.

Orient $C$ and consider its $t$ edges $p_iq_i:\alpha_i$ directed
from $X$ to $Y$. At least $t-2$ satisfy
$up_i\in E(G_c)$ and $vq_i\in E(G_d)$. Deleting such an edge and
attaching $u,v$ gives a strong rainbow spanning path with endpoints
$u,v$ and unused colour $\alpha_i$. Its closure would lift, so
$uv\notin E(F_{\alpha_i})$. Indeed, if
$uv\in E(F_{\alpha_i})$, choose a representative
$f\in H_{\alpha_i}$ for the closing pair $uv$ and prescribe it.
The other $n-1$ pairs of the resulting ordinary Hamilton cycle are
strong, so their codegrees are at least $r-1$. The lifting lemma
(Lemma~\ref{merge:lifting}) then gives a rainbow Berge Hamilton
cycle, a contradiction. Inserting $u$ at $p_iq_i$, using
colour $c$ on $up_i$, also shows that $uq_i$ has codegree at most
one in $H_{\alpha_i}$, by Proposition~\ref{prop:even-one}.

Since $\binom{t-2}2>9t/4$, the bounds on the numbers of nonedges
give two such indices $i,j$ with
$p_ip_j\in E(G_c)$ and $q_iq_j\in E(G_d)$. The chords join $q_i\cdots p_j,q_j\cdots p_i$ into a strong
rainbow $2t$-cycle with $\alpha_i,\alpha_j$ unused. Both unused shadows omit $uv$, so their $A$-sets at $u$
coincide as a parity class. The adjacent vertices $q_i,q_j$ lie
outside this class, a contradiction.

After interchanging $u,v$ if needed, we claim
$|B_c(u)\cap A_c(v)|\ge t-2$.
If both $A$-sets at each omitted vertex lie in one parity class,
the preceding argument makes these classes coincide, proving the claim. Otherwise, choose $u$ whose two $A$-sets are not
both contained in one parity class. Then
$A_c(u)=B_c(u)$ has size $t-1$. Put $T=A_c(u)\cup\{v\}$.
Each of the $t+1$ vertices outside $T\cup\{u\}$ occurs in at most
one edge of the colour-$c$ link of $u$. Hence at most $t$ of the
$k$-subsets of $T$ are missing from that link. If two vertices of
$A_c(u)$ lie outside $A_c(v)$, at most two link edges contain $v$
together with at least one of them. This would leave at least
\[
 \binom{t-2}{k-2}+\binom{t-3}{k-2}-2\ge2t-7>t
\]
missing $k$-sets, a contradiction. Hence the claim.

No edge of $C$ joins $B_c(u)\cap A_c(v)$ to $B_d(u)$. If all
distinct vertices chosen from opposite sets had distance at least
three, Lemma~\ref{lem:separation}(ii) would give
$2t\ge(t-2)+(t-1)+4$ when the sets differ, and
$2t\ge3(t-1)$ when they coincide. Both are impossible.
Thus $C$ has an arc $x,z,x'$ with
$xu\in E(G_c)$, $x\in A_c(v)$ and $x'u\in E(G_d)$.
Replace $z$ by $u$ using colours $c,d$. Let $i$ be the released
colour of $xz$. Then $zv\notin E(G_i)$: otherwise $xv:c,vz:i$ insert $v$
into the original cycle, giving a rainbow ordinary $(n-1)$-cycle
excluded by Proposition~\ref{prop:even-one}.

Now take this new cycle as $C$, relabelling so
$uv\notin E(G_c)$. Both omitted vertices have the same colour-$c$
strong neighbourhood, a parity class $X$ containing both
colour-$d$ $A$-sets.
Let $z$ belong to the other parity class, and let $i,j$ be its two
incident colours. Insertion using colour $c$ excludes both $zu$
and $zv$ from $G_i\cup G_j$. Replace $z$ by $u$ using colours
$c,d$; this is possible because $B_d(u)$ omits at most one of the
two neighbours of $z$. The new omitted pair $zv$ is absent from
$G_i$ and $G_j$. Both new omitted vertices have a parity class as their strong
neighbourhood in colours $i,j$. For $z$ it is $X$, since $zu$ is
absent and $u$ occupies its former place; the class for $v$ is the same.
Vary $z$ and interchange $u,v$: both original omitted vertices
have strong neighbourhood $X$ in every $C$-colour.

For each $z$ in the other parity class, replace $z$ by $u$ once
more, now using the two colours originally incident with $z$.
This leaves $c,d$ unused. Since the original $A$-sets at $u,v$
are contained in $X$, the pairs $zu,zv$ are absent from both unused
strong shadows. The new omitted vertices therefore have strong
neighbourhood $X$ in colours $c,d$ and, by the preceding argument
applied to this cycle, in every colour it uses. Interchanging $u,v$ gives
the same neighbourhoods at $u$. Thus every strong shadow contains
the complete bipartite graph between $X$ and
$V\setminus X$, contrary to Lemma~\ref{lem:shadow-degree}.

\end{proof}

\section{Extending a Berge cycle of length \texorpdfstring{$n-1$}{n-1}}\label{sec:penultimate}
\begin{lemma}\label{lem:released-exclusions}
Let $\mathcal H=(H_1,\ldots,H_n)$ be an $r$-graph system on $n$ vertices,
$3\le r<n$. For $\delta\ge0$, suppose
\begin{equation}\label{eq:section5-assumptions}
 \delta_1(H_i)\ge\delta\ (i\in[n]),\qquad
 \mathcal H\text{ has no rainbow Berge Hamilton cycle}.
\end{equation}
For a maximizing rainbow Berge $(n-1)$-cycle $C$, put $m=n-1$.
Each trace is independent in the ordinary core of $C$ and avoids
the marked ends. Every edge of $H_\alpha\setminus R$ contains at
most one marked tail and at most one marked head.
If $\lambda=2$ and $|U|=m-3$, then
$\Gamma$ is not the disjoint union of two equal cliques. For $a=v_i\in U$,
\begin{equation}\label{eq:released-exclusions}
 T_i\cap(\{w\}\cup N_\Gamma(a)^+\cup Y^+)=\varnothing.
\end{equation}
The reverse orientation gives the corresponding exclusion. Thus
every trace $A$ satisfies
\begin{align}
 v_i\in A&\ \Longrightarrow\ h_i\cap(A^+\setminus\{v_{i+1}\})=\varnothing,
 \label{pen:eq:two-cut-plus}\\
 v_{i+1}\in A&\ \Longrightarrow\ h_i\cap(A^-\setminus\{v_i\})=\varnothing.
 \label{pen:eq:two-cut-minus}
\end{align}
If $U_a$ is the union of traces containing $a$ and $d=d_\Gamma(a)$, then
\begin{equation}\label{eq:released-local-union}
 |U_a|=1+d\le m-r-\lambda+1.
\end{equation}
Moreover, $d_R(a^\pm)\le m-d$ and
$\delta\le\binom{m-d-\lambda-1}{r-1}+m-d-1$.
In particular, $n=2r+2$, $\delta\ge r+1$ and $\lambda=0$ imply
$\Delta(\Gamma)\le r$.
\end{lemma}
\begin{proof}
A trace containing an unmarked core pair would allow its unused
hyperedge to replace the representative of that pair, with the
unused colour assigned to it. The old representative avoids $w$,
whereas the trace edge contains $w$, so this would increase
$d_R(w)$, a contradiction. If a trace meets a marked end $a$,
its edge $f$ and the marked representative $h$ on the pair $ab$
can be assigned to $aw$ and $wb$, respectively, inserting $w$ and
giving a rainbow Berge Hamilton cycle. Thus every trace is independent
in the ordinary core of $C$, and $U$ is disjoint from all marked ends.
By the two-mark exclusion, every unused hyperedge in an unused colour
contains at most one marked tail and at most one marked head.

Let $a=v_i\in U$, $b=a^+$, and let $h_i\in H_{c_i}$ be the
representative assigned to $ab$. Choose
$g\in H_{c_i}(b)\setminus(R\setminus\{h_i\})$ and an unused
colour $\alpha$. Choose $f\in H_\alpha\setminus R$ containing
$w,a$. Trace independence gives $b\notin f$, so $f\ne g$.
If $w\in g$, assign $f$ to $aw$ and $g$ to $wb$; this inserts
$w$ and gives a Hamilton cycle, a contradiction.

Now let $p\in N_\Gamma(a)$, and choose
$f\in H_\alpha\setminus R$ containing $w,a,p$. If $p^+\in g$,
replace the core pairs $ab,pp^+$ by $ap,bp^+$, using $f$ on
$ap$ and $g$ on $bp^+$. The two removed representatives are
unmarked, whereas $f$ contains $w$. The resulting $(n-1)$-cycle
therefore has one more representative containing $w$, contrary
to maximality. If $g$ contains the head $v$ of a marked pair $uv$,
with representative $h\in H_\gamma$, replace $ab,uv$ by $aw,wu,bv$,
using $f,h,g$, respectively. This gives a rainbow Berge Hamilton
cycle; the same construction works when $b=u$.
These arguments prove \eqref{eq:released-exclusions}. Taking
$g=h_i$ and reversing the orientation gives
\eqref{pen:eq:two-cut-plus} and \eqref{pen:eq:two-cut-minus}.

The successor map is injective, and none of the $d$ forbidden
successors is a marked head: its predecessor lies in $U$, which
avoids marked ends. They also avoid $b$, since $N_\Gamma(a)$ has
no loop. Thus the $d$ forbidden successors, the $\lambda$ marked
heads and $b$ are pairwise disjoint. Hence
$r-1\le|T_i|\le m-1-d-\lambda$, proving
\eqref{eq:released-local-union}. By \eqref{pen:eq:two-cut-plus},
the $d$ representatives leaving $N_\Gamma(a)$ avoid $b$, so
$d_R(b)\le m-d$. At most $m-d-1$ selected edges other than $h_i$
contain $b$; all other colour-$c_i$ edges there lie in
$T_i\cup\{b\}$. This gives the degree bound; reverse for $a^-$.
For $m=2r+1$ and $\lambda=0$, these bounds give $d\le r+1$,
with equality forcing $\delta\le r$. Hence $\delta\ge r+1$ implies
$d\le r$.

Suppose $\lambda=2$, $|U|=m-3$, and $\Gamma$ is the disjoint
union of two equal cliques. The two marked core pairs form an
ordinary path on the three inactive vertices. Label the remaining
core segment in order as $(1,\ldots,m-3)$. Since each clique is
independent in the core, the two cliques alternate along this
segment. Applied to the representative of $12$, the forward
exclusion forbids the even vertices from $4$ onwards and the marked
heads. The reverse exclusion forbids the odd vertices from $3$
onwards and the marked tails. The representative also avoids $w$,
so it would be contained in $\{1,2\}$, contrary to $r\ge3$.
\end{proof}

\begin{lemma}\label{lem:independent-clique}
Let $\mathcal H=(H_1,\ldots,H_n)$ be an $r$-graph system on
$n=2t+2$ vertices, $3\le r\le t+1$, with
$\delta_1(H_i)>\binom t{r-1}$ for every $i$ and no rainbow Berge
Hamilton cycle. For a maximizing rainbow Berge $(n-1)$-cycle $C$,
if $\lambda>0$, then $K_t\nsubseteq\Gamma$. Moreover, its trace
union $U$ cannot be an independent $t$-set in the ordinary core.
\end{lemma}
\begin{proof}
Suppose $\lambda>0$ and $T$ induces $K_t$ in $\Gamma$.
By Lemma~\ref{lem:released-exclusions}, $T$ is independent in
the core, so its complement $S$
spans a single core pair $uv$, the unique mark. Let $h$ be its
representative.

For $a\in T$ and $b=a^+$, the forbidden heads
$(T\setminus\{a\})^+\cup\{v\}$ exhaust $S\setminus\{b\}$.
Thus every edge available at $b$ after releasing $ab$, including
its representative, lies in $T\cup\{b\}$. Reverse to obtain this
restriction for every pair meeting $T$.
Choose $b\in S\setminus h$, since $|h\cap S|\le r-1<t+1$.
Only its two incident representatives can contain $b$. In the
incoming colour, every edge at $b$ therefore lies in $T\cup\{b\}$,
giving degree at most $\binom t{r-1}$.

If $U$ were itself an independent $t$-set in the core, at most one
core pair could be marked. Since the unused colour has degree greater
than $\binom t{r-1}$ at $w$, there must be exactly one marked
representative and every $(r-1)$-subset of $U$ must be a trace.
Thus $U$ induces the forbidden $K_t$ in $\Gamma$.
\end{proof}

\begin{lemma}\label{lem:protected-marked-cut}
Let $\mathcal H=(H_1,\ldots,H_n)$ be an $r$-graph system on
$n=2t+1+\varepsilon$ vertices, $\varepsilon\in\{0,1\}$, with
$3\le r\le t$ and $\delta_1(H_i)\ge D(n,r)$ for every $i$.
Assume $\mathcal H$ has no rainbow Berge Hamilton cycle, and let
$C$ be a maximizing rainbow Berge $(n-1)$-cycle.
For $ap\in E(\Gamma)$, release the representatives of $aa^+,pp^+$
in their colours, giving neighbourhoods $T_a,T_p$ at $a^+,p^+$.
Then $|T_a|+|T_p|\le n-1$. Also $|T_a|\le t$ if $r=3$,
$r<t$, or the forward arc from $p^+$ to $a$ contains no marked pair.
\end{lemma}
\begin{proof}
Put $k=r-1$. Release the representatives of $aa^+,pp^+$ in
colours $i,j$, and join $p$ to $a$ with an unused edge $h\supseteq\{w,p,a\}$.
The resulting path $P$ has core order $a^+Cp,aC^{-1}p^+$,
ends $x=a^+,z=p^+$, and $\lambda+1$ marks at $w$.
Its neighbourhoods $A,B$ in colours $i,j$ are $T_a,T_p$:
$h$ avoids both ends, and an opposite released representative
would close $P$. Lemma~\ref{lem:common-endpoint-cuts}(ii) gives
the sum bound. Hereafter $R$ denotes the current path's representatives.

Suppose $|A|>t$. For any such path retaining $A$ and at least
$\lambda+1$ marks, Lemma~\ref{lem:fixed-end-head-bound} gives
selected degree at most $n-2-|A|\le t-1$ at the opposite end.
Its other unused colour has at least $D(n,r)-t+1>\binom{t-1}k$ unused
incident edges, hence at least $t$ neighbours. Lemma~\ref{lem:common-endpoint-cuts}(ii)
forces $n=2t+2$, $|A|=t+1$, $|B|=t$, with every core pair counted.
Retaining the original colour-$i$ edges at $x$ therefore keeps $A$
unchanged.

Let a selected $g\in H_j(z)$ represent $u\to v$ in colour $c$.
Rotate with $g:j$ on $uz$ to make $v$ the opposite endpoint;
if $v=z$, only recolour. Direct closure excludes $v\in A$.
The equality in the endpoint count at $u\to z$ then gives
$u\in N_c^R(v)$. Return to the original core order and replace $g$ by an unused
colour-$c$ edge on $uv$. Since $v\notin A$, every original
colour-$i$ edge at $x$ remains unused.

For triples, repeat for every nonterminal representative in
$H_j(z)$. Each is $uvz$ and is replaced by an edge avoiding $z$,
so their number decreases. Since $uvz$ avoids $w$, each step
retains $h$ and all marks; it also keeps every previously unused
colour-$j$ edge at $z$ available. The endpoint count keeps $A,B$
fixed. The process
ends with only the terminal representative possibly in $H_j(z)$.
Degree $D(n,r)=\binom t2+1$ forces it outside $B\cup\{z\}$ and makes
all triples containing $z$ and two vertices of $B$ available in colour
$j$. Each neighbour in $B$ has $t-1\ge2$ such edges.
For $3\le k\le t-2$, there are already at least two edges per
neighbour, since $D(n,r)-(t-1)-\binom{t-1}k\ge2$.

In each case degree supplies a selected $g\in H_j(z)$ leaving
$B\cup\{z\}$; perform the same exchange. If $w\in g$, the order
$xPp,w,zP^{-1}a,x$ is Hamilton, using $h$, the released $g:j$,
and an original colour-$i$ edge on $ax$. Otherwise the marks survive.
Releasing $g$ enlarges $B$, so its replacement must remove
an old neighbour. This is impossible in the two-edge cases above.

Finally let $r=t>3$ and suppose the reversed arc is unmarked.
By \eqref{eq:released-exclusions}, $p\to a$ is the only marked pair
with head in $A$. Before and after the exchange, no unused
colour-$j$ edge contains $wz$: together with $h$ and an original
colour-$i$ edge on $ax$, it would be Hamilton. The attachments are distinct: a common edge would contain $xz$
and close $P$. The released $g$ also differs from the original
unused edge.
The changed pair has head outside $A$, so $p\to a$ remains the
unique marked pair with head in $A$. The $t$ tails determined by the opposite endpoint
neighbourhood and $t+1$ heads in $A$ cover all $2t$ core pairs, with a unique
common pair, which must be marked by Lemma~\ref{lem:common-endpoint-cuts}(ii).
Thus $B$ equals $p$ together with the tails whose heads lie outside
$A$, both before and after the exchange, contradicting the new
neighbour supplied by $g$.
\end{proof}

\begin{lemma}\label{turn21:unique-mark-boundary}
Let $\mathcal H=(H_1,\ldots,H_n)$ be an $r$-graph system on
$n=2t+1+\varepsilon$ vertices, $\varepsilon\in\{0,1\}$, with
$3\le r\le t$ and $\delta_1(H_i)\ge D(n,r)$ for every $i$.
Assume $\mathcal H$ has no rainbow Berge Hamilton cycle, and let
$C$ be a maximizing rainbow Berge $(n-1)$-cycle.
If $\lambda=1$, both marked ends and their two outer core neighbours
lie outside $U$. Hence at least $D(n,r)-1$ traces lie on at most
$n-5$ vertices.
\end{lemma}
\begin{proof}
Suppose an active vertex $a$ immediately precedes the marked
core pair $bc$, whose representative is $h\in H_\gamma$.
Let $g\in H_i$ be the representative assigned to $ab$, and choose
an unused colour $\alpha$ and an edge
$f\in H_\alpha\setminus R$ containing $w,a,p$ for some
$p\notin\{a,w\}$. Replace the core segment $a,b,c$ by $a,w,c$,
using $f$ on $aw$ and $h$ on $wc$. The new cycle omits $b$.
Since $\lambda=1$ and $h$ contains $b$, no other representative
of this cycle can contain $b$; otherwise its omitted vertex $b$
would have selected degree at least two, contrary to the global
maximality of $\lambda$.

In the released colour $i$, at most one selected representative
of the new cycle contains $b$, namely $h$. Hence at least
$D(n,r)-1\ge\binom t{r-1}$ unused hyperedges of $H_i$ contain
$b$. If their traces at $b$ had union of size at most $t-1$, there
would be at most $\binom{t-1}{r-1}$ such hyperedges. Thus the trace
union $A$ at $b$ has size at least $t$. If $|A|=t$, all
$\binom t{r-1}$ possible traces occur, so every $(r-1)$-subset
of $A$ is a trace. Since every trace is independent in the new
core, $A$ is independent. Lemma~\ref{lem:independent-clique}
excludes this in even order. In odd order, the core has $2t$ vertices,
so an independent set of size $t$ is a maximum independent set of
$C_{2t}$. Its complement is the other alternating class and is also
independent, but contains the core pair $wc$, a contradiction.
Therefore $|A|>t$. The released edge $g$ gives $a\in A$,
whereas the marked ends $w,c$ do not belong to $A$.

Releasing $ab$ on the original cycle gives the same unused
$H_i$-neighbourhood $A$ at $b$, since the trace of $f$ is independent
and therefore avoids the adjacent vertex $b$. In the cycle order
$a,b,c,\ldots,a$, the trace vertex $p$ avoids the marked ends
$b,c$ and is nonadjacent to $a$; hence it lies on the open arc
from $c$ to $a$. The forward arc from $p^+$ to $a$ therefore
avoids $bc$. Lemma~\ref{lem:protected-marked-cut} gives
$|A|\le t$, a contradiction. Reversing the orientation excludes
the other outer neighbour.
\end{proof}

\begin{lemma}\label{lem:zero-neighbour-load}
Let $\mathcal H=(H_1,\ldots,H_n)$ be an $r$-graph system on
$n=2t+1+\varepsilon$ vertices, $\varepsilon\in\{0,1\}$, where
$1\le k=r-1\le t$. Suppose $\delta_1(H_i)>\binom tk$ for every $i$,
no rainbow Berge Hamilton cycle exists, and every rainbow Berge
$(n-1)$-cycle has no representative containing its omitted vertex.
For any such cycle $C$ with representatives $R$ and trace union $U$,
$d_R(v)\le t+\varepsilon$ for $v\in N_C(U)$. If $r\ge3$, this
improves to $t+\varepsilon-1$. If also $r\le t$, both core neighbours
of $v$ lie in $U$ and $t+\varepsilon\ge4$, it improves to
$t+\varepsilon-2$.
\end{lemma}
\begin{proof}
Put $m=n-1=2t+\varepsilon$. Label the core vertices in cyclic
order so that $v_1\in U$, and write
\[
 C=v_0h_0v_1h_1\cdots v_{m-1}h_{m-1}v_0,
\]
with indices modulo $m$. Delete the core pair $v_0v_1$ and its
representative $h_0$. Let $\alpha$ be the unused colour and
$\gamma$ the colour of $h_0$. Every edge containing $w$ is unused,
since $\lambda=0$. On the resulting rainbow Berge path with core
order $(v_1,\ldots,v_{m-1},v_0)$, let $S$ be the set of heads of
core pairs whose representatives contain $v_0$. Then $v_0\in S$
and $|S|=d_R(v_0)-1$.
If $|S|=1$, put $K=S$. Otherwise, put
$K=S\cup\{v_{i+1}:v_i\in S\setminus\{v_0\}\}
\cup\{v_{m-1}\}$.
In both cases, $v_1\notin K$.

For every $v_j\in K$, edges containing $w,v_1$ and $w,v_j$ in
opposite colours $\alpha,\gamma$ must coincide. If distinct, they
close a rainbow Berge Hamilton cycle directly when $v_j=v_0$.
For $v_j\in S\setminus\{v_0\}$, use the core order
$w,v_1,\ldots,v_{j-1},v_0,v_{m-1},\ldots,v_j,w$, assigning the
former representative $h_{j-1}$ of $v_{j-1}v_j$ to $v_{j-1}v_0$.
Using an attachment at $v_{j+1}$ instead omits $v_j$ while
retaining that representative containing $v_j$, contrary to
$\lambda=0$. Finally, when $|S|\ge2$, the order
$w,v_1,\ldots,v_{m-1},w$ omits $v_0$ and retains a representative
containing $v_0$, giving the same contradiction. These arguments
apply with the two colours interchanged.

Fix a colour-$\alpha$ edge $f$ containing $w,v_1$. If
$|K|\ge t+\varepsilon$, the degree bound supplies a colour-$\gamma$
trace meeting $K$. Its edge is $f$; every other colour-$\alpha$
trace avoids $K\cup\{v_1\}$, so there are at most
$\binom{t-1}k<(\binom tk+1)-1$ such traces. This contradicts the degree
at $w$. Hence $|K|\le t+\varepsilon-1$, and therefore
$d_R(v_0)\le t+\varepsilon$.

A stronger bound holds if a colour-$\alpha$ trace meets $K$ and
avoids $v_1$: no colour-$\gamma$ trace can contain $v_1$. The
fixed edge $f$ also excludes colour-$\gamma$ traces meeting $K$.
Thus at least $\binom tk+1$ colour-$\gamma$ traces avoid
$K\cup\{v_1\}$, forcing $|K|\le t+\varepsilon-2$.
Now assume $3\le r\le t+1$. If
$d_R(v_0)=t+\varepsilon$, then
$|S|=t+\varepsilon-1$ and the preceding bound forces $K=S$.
Thus $S$ is a terminal interval in the core order, closed under
successors up to $v_0$. A trace containing $v_1$ avoids $v_0$ by
independence. It also avoids every $v_j\in S\setminus\{v_0\}$ by
\eqref{pen:eq:two-cut-minus}, because the representative entering
$v_j$ contains $v_0=v_1^-$. The stronger bound excludes any other
trace meeting $K$. Deleting $K$ leaves an ordinary core path on
$t+1$ vertices, with only
$\binom{t-k+2}k\le\binom tk$ independent $k$-sets, fewer than
the required $\binom tk+1$ traces. Thus
$d_R(v_0)\le t+\varepsilon-1$.

Finally, suppose $v_{m-1},v_1\in U$, $r\le t$,
$t+\varepsilon\ge4$ and
$d_R(v_0)=t+\varepsilon-1$. Write $h_j$ for the representative
of the core pair $v_jv_{j+1}$. If $h_{m-2}$ avoids $v_0$, then
$K\setminus S=\{v_{m-1}\}$. Closure under successors up to
$v_{m-1}$ forces the $t+\varepsilon-3$ nonincident core pairs
whose representatives contain $v_0$ to have indices
$m-(t+\varepsilon)+1,\ldots,m-3$. In particular, $h_1$ avoids
$v_0$. Reversing the orientation places those same indices in
the disjoint interval $2,\ldots,t+\varepsilon-2$, a contradiction.
Therefore $h_{m-2}$ contains $v_0$; by reversal, so does $h_1$.
The latter excludes the trace pair $v_{m-1}v_1$ by
\eqref{pen:eq:two-cut-plus}. A trace containing $v_{m-1}$ therefore
meets $K$ and avoids $v_1$, forcing
$|K|\le t+\varepsilon-2$. Hence $K=S$ is terminal. This is
impossible because $v_2\in S$, completing the proof.
\end{proof}

\subsection{Traces on the ordinary core cycle}

\begin{lemma}\label{new:even-distance-two}
Let $\mathcal H=(H_1,\ldots,H_n)$ be an $r$-graph system on
$n=2t+1+\varepsilon$ vertices, $\varepsilon\in\{0,1\}$, with
$3\le r\le t$ and $\delta_1(H_i)\ge D(n,r)$ for every $i$.
Assume $\mathcal H$ has no rainbow Berge Hamilton cycle, and let
$C$ be a maximizing rainbow Berge $(n-1)$-cycle.
Suppose $\lambda=0$. Then
\begin{enumerate}
\item[(i)] No trace contains two vertices at core distance two.
\item[(ii)] If $2k\ge t+1+\varepsilon$ or $(n,r)=(8,3)$,
then no such cycle $C$ exists. This includes $r=t$.
\end{enumerate}
\end{lemma}
\begin{proof}
For (i), suppose an unused colour $\alpha$ has a hyperedge
$f\in H_\alpha\setminus R$ containing $w$ and the ends
$u,v$ of a core segment $u,z,v$. Replace that segment by the pair
$uv$, represented by $f$, obtaining a rainbow Berge $(n-2)$-cycle
$C'$. Write
$R'=E(C')$, and let $i,j$ be the released colours of the former
representatives of $uz,zv$.
Neither colour has an edge containing $w,z$: together with $f$
it would insert $w$ at the corresponding original core pair.
The edges are distinct, since $z\notin f$ and every original
representative avoids $w$.

Delete a core pair $ab$ of $C'$ whose representative $h$ contains
exactly one of $w,z$, and let $q$ be the other vertex. This produces
a rainbow Berge path $P$. Reserve $h$ for closure and keep $R'$ as the exclusion set in
colours $i,j$, also when opening $C'$ at another pair. If $t=3$, we initially open only at $uv$, so $h=f$
and $q=z$.

Put $B=\binom tk$ and $L=B-t+4$. In each of colours $i,j$, at
least $L$ edges containing $q$ lie outside $R'$. For $q=w$, only
$f$ among the selected representatives contains $w$. For $q=z$,
Lemma~\ref{lem:zero-neighbour-load} bounds its original selected
degree by $t+\varepsilon-2$ when $t\ge4$, and by two when $t=3$;
both original representatives incident with $z$ have been removed.

Put $A_i=N_i^{R'}(q)$ and $A_j=N_j^{R'}(q)$.
For consecutive $x\in A_i,y\in A_j$, distinct attachments insert
$q$ and close with $h$, contradicting $\lambda=0$. Hence both
adjacencies have the same unique representative in both colours.

In a maximum matching of the ordinary core of $P$, let $d$ count
pairs with opposite-colour ends. By uniqueness, their $2d$ ends
belong to both neighbourhoods and have trace degree one in each colour. Each other matched pair, and the possible
unmatched vertex, contributes at most two memberships. Since
$|V(P)|=2t-1+\varepsilon$, one colour has at most $t-d$ neighbours
outside those pairs and at most $d$ traces meeting them.
Put $s=t-k$. If $1\le d\le s$, then $L\le d+\binom{t-d}k$.
By convexity the right side is at most
$\max\{1+\binom{t-1}k,s+1\}<L$, since
$L\ge\binom{t-1}k+3\ge s+3$.
If $d>s$, every trace needs at least $d-s$ of the $2d$ ends, giving
$L\le\min\{d,2d/(d-s)\}\le s+2$.
Hence $d=0$. The degree bound gives $|A_i|,|A_j|\ge t$, and the
matching count forces equality. This excludes $t=3$, since
then $L=4>B$. Every neighbour has at least
$L-\binom{t-1}k\ge3$ distinct available representatives. No vertex of $A_i$ is consecutive to a vertex of $A_j$
in the core path. Lemma~\ref{lem:separation}(iii) therefore gives
$A_i=A_j=A$, where $A$ is an independent $t$-set.

Degree supplies a colour-$i$ edge containing $q$ outside $A\cup\{q\}$. It represents a core pair $xy$, since all available edges lie inside
$A\cup\{q\}$ and $h$ avoids $q$.
If $x\in A$, use this representative on $qy$, an available
colour-$i$ edge on $xq$, and $h$ to close an $(n-1)$-cycle with
$h$ containing its omitted vertex. Thus $xy$ is disjoint from $A$.
The independent set $A$ meets $2t-|A\cap\{a,b\}|$ of the
$2t-2+\varepsilon$ core pairs of $P$. Because $A$ is independent, no path edge has both ends in $A$.
Each internal vertex of $A$ is incident with two path edges, while
an endpoint of $P$ in $A$ is incident with one. Thus the number of
core edges of $P$ meeting $A$ is
$2|A|-|A\cap\{a,b\}|=2t-|A\cap\{a,b\}|$.
Since $P$ has $2t-2+\varepsilon$ core edges and the pair $xy$
is disjoint from $A$, we obtain
\[
1\le(2t-2+\varepsilon)-(2t-|A\cap\{a,b\}|)
=\varepsilon-2+|A\cap\{a,b\}|.
\]
As $\varepsilon\in\{0,1\}$ and $|A\cap\{a,b\}|\le2$,
this forces $\varepsilon=1$ and $a,b\in A$.

Now open $C'$ first at $uv$, with $q=z$, and choose a core pair
$xy$ disjoint from the resulting set $A$. Write the core order as
$uP_1x,yP_2v$. The representative of $xy$ contains $z$ and avoids
$w$. Next open $C'$ at $xy$, with $q=w$, retaining the exclusion
set $R'$ and the unused colours $i,j$.
The argument gives $v\in N_i^{R'}(z)$ and
$y\in N_j^{R'}(w)$. Use the representative of $xy$ on $xz$, the
edge $f$ on $wu$, and available edges on $zv$ and $yw$ in colours
$i$ and $j$, respectively. These form a rainbow Berge Hamilton
cycle with core order $uP_1x,z,vP_2^{-1}y,w,u$.
The available edges differ, since neither colour contains an edge
on $wz$. This proves (i).

For (ii), the cyclic gaps of a trace are at least three, so no
trace exists if $3k>m$. Otherwise, distinguish one trace vertex
and distribute $m-3k$ among its $k$ gaps; each trace is counted $k$ times. When $2k\ge t+1+\varepsilon$,
the number of traces is at most
\[
 \frac{m}{k}\binom{m-2k-1}{k-1}\le
 \frac{m}{k}\binom{t-2}{k-1}
 =\frac{m(t-k)}{t(t-1)}\binom tk\le B
\]
fewer than the $D(n,r)=B+1$ traces required by the degree
at $w$.

It remains to consider $(n,r)=(8,3)$. The graph on $V(C)$ joining vertices at core distance three is
an ordinary seven-cycle. Every trace is one of its edges, so the
required four traces give $|U|\ge5$. If $|U|=5$, the two inactive vertices must
be adjacent in the auxiliary seven-cycle, since deleting two
nonadjacent vertices leaves only three edges. Consequently no core
vertex has two inactive core neighbours, also when $|U|>5$.
Lemma~\ref{lem:zero-neighbour-load} bounds every selected degree
by three, and bounds it by two at each of the $2|U|-7$ vertices
with two active core neighbours. Thus the selected degrees sum to less than $21$. All seven
representative triples avoid the omitted vertex, so the sum is
$21$, a contradiction.
\end{proof}

\begin{lemma}\label{new:near-complete-triangle-capacity}
If $t\ge5$,
$3\le k<t$, $0\le q\le t-2$ and $G$ is a $K_t$-free graph on at most
$t+q$ vertices with $\Delta(G)\le t-1$, then its number of $k$-cliques satisfies
\begin{equation}\label{eq:neighbourhood-clique-bound}
 c_k(G)\le\max\left\{\binom{t-1}k+\binom{q+1}k,\;
 \binom tk-\binom{t-3}{k-3}-\binom{t-4}{k-3}+\binom qk\right\}.
\end{equation}
\end{lemma}
\begin{proof}
Let $a(v)$ count the $k$-cliques avoiding $v$ but meeting $N(v)$.
Chao and Dong~\cite[(3),(4)]{CD} give
$\sum_v a(v)\le\sum_v\binom{d(v)}k$. Choose $v$ with $a(v)\le\binom{d(v)}k$. Deleting $N[v]$ removes
at most $\binom{d(v)}k+c_{k-1}(G[N(v)])$ cliques. Repeat in
induced subgraphs: a step deleting $b$ vertices removes at most
$\binom bk$ cliques. If every step has $b\le t-1$, at most
$t+q\le2(t-1)$ vertices are deleted, so convexity gives the first bound.

Otherwise, consider a step with $|N[v]|=t$. Since the current
graph is $K_t$-free, $N(v)$ has a missing pair.
Two missing pairs exclude at least
$\binom{t-3}{k-3}+\binom{t-4}{k-3}$ of its $(k-1)$-sets:
count sets containing the first pair, then those containing the
second and avoiding an endpoint of the first outside it.
This step removes at most
$\binom tk-\binom{t-3}{k-3}-\binom{t-4}{k-3}$ cliques.

If $xy$ is the unique missing pair, every other vertex of $N[v]$
already has degree $t-1$, and $x,y$ have at most one outside
neighbour each. No triangle, hence no $k$-clique, crosses $N[v]$. Exactly
$\binom tk-\binom{t-2}{k-2}$ cliques are removed, within the
same bound by Pascal's identity. All other steps, including earlier
ones, delete at most $q$ vertices in total. Their aggregate
contribution is at most $\binom qk$ by convexity, proving the second bound.
\end{proof}

\begin{lemma}\label{merge:uniform-zero}
Let $\mathcal H=(H_1,\ldots,H_n)$ be an $r$-graph system on
$n=2t+1+\varepsilon$ vertices, $\varepsilon\in\{0,1\}$, with
$4\le r<t$, $\delta_1(H_i)\ge D(n,r)$ for every $i$, and no
rainbow Berge Hamilton cycle. Put $k=r-1$. Every maximizing rainbow
Berge $(n-1)$-cycle satisfies
\begin{equation}\label{eq:uniform-released-degree}
 \Delta(\Gamma)+\lambda\le t+\varepsilon-1.
\end{equation}
If $2k\ge t$, every rainbow Berge $(n-1)$-cycle has no representative
containing its omitted vertex. If $2k\ge t+1+\varepsilon$, no such
cycle exists.
\end{lemma}
\begin{proof}
Write $m=2t+\varepsilon$.
Since $D(n,r)>m$, a trace exists unless $(n,r)=(12,4)$.
An empty trace family there marks all eleven pairs.
There are at least $12D(n,r)/4-m=22$ unused colour-$\alpha$ edges,
each avoiding $w$ and thus containing four marked tails. This
contradicts Lemma~\ref{lem:released-exclusions}. Hence a trace
exists in every case.

Put $d=\Delta(\Gamma)\ge k-1$. If $d+\lambda\ge t+\varepsilon$,
Lemma~\ref{lem:released-exclusions} gives
$D(n,r)\le\binom{t-1}k+m-d-1$. For $t\ge6$, this contradicts
$\binom{t-1}{k-1}\ge\binom{t-1}2\ge m-k,$
since $d\ge k-1$. When $t=5$, the same inequality gives
$d\le2+\varepsilon$ and $\lambda\ge3$.
A graph of maximum degree at most three has at most as many triangles
as vertices, by counting triangle--vertex incidences. Also, the marked
core pairs have at least $\lambda+1$ distinct inactive ends.
It follows that
\[
 11-\lambda\le |\mathcal T|\le c_3(\Gamma)
 \le |U|\le10-\lambda,
\]
a contradiction. This proves \eqref{eq:uniform-released-degree}.

Now suppose $\lambda>0$ and $2k\ge t$.
First consider $\lambda=1+\varepsilon+j$ for an integer $j\ge0$.
Then $|U|\le m-\lambda-1\le2(t-1)$ and
$k-1\le\Delta(\Gamma)\le t-2-j$. Add isolated vertices so that
$\Gamma$ has exactly $2(t-1)$ vertices. Count incidences of vertices with $k$-cliques. Decreasing the upper
argument of $\binom{t-2}{k-1}$ by $j$ reduces the coefficient by
at least $j$, giving
\[
 \binom tk-\varepsilon-j\le c_k(\Gamma)
 \le 2\binom{t-1}k-\frac{2(t-1)}k j
 \le \binom tk-\frac{2(t-1)}k j .
\]
Since $2(t-1)/k>2$, these inequalities force $j=0$.
They also force $2k=t$: otherwise
$\binom tk-2\binom{t-1}k=(2k-t)\binom tk/t>1$.
When $2k=t$, the total deficit from the maximum
$\binom{t-2}{k-1}$ $k$-cliques at each vertex is at most $k$.
A zero-deficit vertex lies in a $K_{t-1}$ component: $k\ge3$ and
$\Delta(\Gamma)\le t-2$ force its neighbourhood to be a clique.
The vertices outside these components number a multiple of $t-1>k$
and each contributes at least one to the deficit. Hence none remains,
and $\Gamma=2K_{t-1}$. In odd order this contradicts
Lemma~\ref{turn21:unique-mark-boundary}; in even order it
contradicts Lemma~\ref{lem:released-exclusions}.

It remains to consider even order with $\lambda=1$.
By Lemma~\ref{turn21:unique-mark-boundary}, the trace family has at least $\binom tk$ members
on at most $2t-3$ vertices. Also $\Delta(\Gamma)\le t-1$,
and Lemma~\ref{lem:independent-clique} implies that $\Gamma$ is
$K_t$-free. Apply Lemma~\ref{new:near-complete-triangle-capacity}
with $q=t-3$. Both bounds are below $\binom tk$,
since $2k\ge t$ gives
\[
 \binom{t-1}k+\binom{t-2}k<\binom tk,
 \qquad
 \binom{t-3}k\le\binom{t-3}{k-3}.
\]
Thus $\lambda=0$, so every $(n-1)$-cycle has no representative
containing its omitted vertex. Lemma~\ref{new:even-distance-two}(ii)
gives the final assertion.
\end{proof}

\section{Systems with \texorpdfstring{$n\le2r+2$}{n at most 2r+2}}\label{sec:large}
\begin{lemma}\label{rt:lem:orbit}
Let $\mathcal H=(H_1,\ldots,H_n)$ be an $r$-graph system on $n$ vertices,
$3\le r<n$. Fix one endpoint of a rainbow Berge path of order at least two. Let $R$ be its
representative set, and let $X$ consist of the opposite endpoints obtainable
by rotations preserving $R$ and its colour assignment. Then
$d_R(x)\le |X|$ for every $x\in X$. If equality holds for every
$x\in X$, then $|X|<r$. Moreover, if $n\le2r$ and
$\delta_1(H_i)\ge r$ for every $i\in[n]$, then, for every unused
colour $a$, some edge in $H_a\setminus R$ meets $X$.
\end{lemma}
\begin{proof}
Orient an attainable path from $x\in X$ towards the fixed end.
The representatives containing $x$ permit rotations to their distinct
tails, giving $d_R(x)\le|X|$.
If equality holds throughout $X$, these tails exhaust $X$. Let $z$
be the last one; then $z^+\notin X$. Rotations to earlier tails
preserve $zz^+$, and equality forces its representative to contain
every new end. It contains $X\cup\{z^+\}$, so $|X|<r$.

If all unused colour-$a$ edges avoid $X$, then $d_R(x)\ge r$
on $X$, so the first assertion and its equality case give
$|X|\ge r+1$. For $n\le2r$, fewer than $r$ vertices lie outside
$X$, forcing $H_a\setminus R=\varnothing$. This contradicts
$|H_a|\ge n>|R|$, obtained by summing degrees.
\end{proof}

\begin{lemma}\label{lem:near-large}
Let $\mathcal H=(H_1,\ldots,H_n)$ be an $r$-graph system on $n$ vertices,
$3\le r<n\le2r+2$, with $\delta_1(H_i)\ge D(n,r)$ for every $i$.
Assume no rainbow Berge Hamilton cycle exists. If a rainbow Berge
$(n-1)$-cycle exists, then $n\ge2r$. Moreover, if $\lambda$
denotes the maximum of $d_R(w)$ over all such cycles $C$, with
representative set $R$ and omitted vertex $w$, then $\lambda=0$
unless $(n,r)=(8,3)$; at $(n,r)=(8,3)$,
$\lambda\in\{0,2,3\}$.
\end{lemma}
\begin{proof}
Let $C$ be a maximizing Berge $(n-1)$-cycle, with unused colour $\alpha$.
By Lemma~\ref{lem:released-exclusions}, a trace at $w$ is independent
and avoids marked ends. Its $r-1$ vertices meet $2r-2$ unmarked
pairs, so $\lambda\le n-2r+1$ whenever a trace exists.

Suppose $H_\alpha(w)\setminus R=\varnothing$. Then
$\lambda\ge D(n,r)>n-2r+1$, so every maximizing cycle has empty trace family.
Lemma~\ref{lem:released-exclusions} forces each unused edge to meet
at most one marked tail and one marked head.
There are at least $\lceil nD(n,r)/r\rceil-m$ unused edges:
at least one in all cases, and four if $n=2r+2$.
These edges avoid $w$, so the marks cannot cover $C$.
The marks form $c\ge1$ maximal proper paths. The set $Z$ outside their
ends has size $m-\lambda-c\le r-1$; every unused edge needs
at least $r-2$ vertices of $Z$.

If $|Z|=r-2$, every unused edge consists of $Z$, the first vertex
of a marked run and the last vertex of a marked run. There are $c^2$
candidates. Since $\lambda\ge D(n,r)$, we have $c\le1$ for $n\le2r+1$
and $c\le2$ otherwise, so the lower bound forces all candidates
to occur. They cover every unmarked pair along the segment through $Z$.
Assign one candidate hyperedge $f$ to such a pair, preserving the marks.
Let $i$ be the colour released by the former representative. In colour
$i$ the degree-sum bound requires at least $c^2$ unused edges, while
the same two-mark classification allows only the $c^2$ candidate
hyperedges. Hence all candidates would have to be unused in colour $i$.
But $f$ is now selected; whether or not $f$ belongs to $H_i$, at most
$c^2-1$ candidate hyperedges can be unused there. This is a contradiction.

If $|Z|=r-1$, then $n=2r+2$, $\lambda=r+1$ and $c=1$.
Let $S$ be the vertex set of the unique unmarked segment, with
interior $Z$, and let $J=V(C)\setminus S$, so $|J|=r$.
Every vertex of $J$ is internal on the unique marked path and is
therefore both a marked tail and a marked head. The two endpoints
of $S$ are the two boundary marked ends. By the two-mark exclusion,
an unused edge cannot contain a vertex of $J$ together with an
endpoint of $S$, and cannot contain two vertices of $J$. Since every
unused edge needs at least $r-2$ vertices of $Z$, it follows that
each unused edge is either an $r$-subset of $S$ or has the form
$F_u=Z\cup\{u\}$ for a unique $u\in J$.

Every vertex of $J$ occurs in at most one unused edge, so every
vertex occurring at least twice in the unused family lies in $S$.
If all $r+1$ vertices of $S$ occur at least twice, they determine
$S$. Otherwise, at most $r$ vertices of $S$ occur at least twice.
There can then be at most two distinct $r$-subsets of $S$ in the
unused family: any three distinct $r$-subsets of the $(r+1)$-set
$S$ omit three different vertices, so every vertex of $S$ would
occur in at least two of them, a contradiction. Since the unused
family has at least four edges, at least two distinct edges $F_u,F_v$
occur. Their intersection is exactly $Z$. More intrinsically, $Z$ is the
set of vertices occurring in at least three unused edges: each vertex
of $Z$ is absent from at most one unused edge, whereas each endpoint
of $S$ occurs in at most two and each vertex of $J$ in at most one.
The two core neighbours of $Z$ recover the endpoints of the same
unmarked segment $S$.

The $\lambda=r+1$ marked representatives have at most
$(r+1)(r-1)<r^2$ incidences with $J$. Choose $u\in J$ contained
in at most $r-1$ of them. Since $u$ is internal on the marked path,
its two incident core pairs are marked; reassign their representatives
to the two pairs obtained by replacing $u$ with $w$. This gives a
rainbow Berge $(n-1)$-cycle omitting $u$ and retaining the same
representative set $R$.

If $F_u$ is unused in colour $\alpha$, then among the at least
$r+1$ hyperedges of $H_\alpha$ through $u$, all but $F_u$ must be
selected representatives: by the classification above, $F_u$ is
the only possible unused edge through $u$. At most $r-1$ of these
selected representatives are marked, by the choice of $u$; hence
one is an unmarked representative on a core pair $xy$ of $S$.
That pair meets $Z$. Use $F_u$ on the pair joining $u$ to its endpoint
in $Z$, and use the selected representative on the other pair of
$x-u-y$. This inserts $u$ into the new cycle and gives a rainbow
Berge Hamilton cycle, a contradiction. If $F_u$ is not unused,
there is no unused edge of colour $\alpha$ through $u$; hence all
at least $r+1=\lambda$ incident hyperedges of $H_\alpha$ are
selected, so $d_R(u)\ge\lambda$. The new cycle is therefore
maximizing.

A maximizing cycle with two marked paths would have
$|Z|=m-\lambda-2=r-2$, which is excluded by the preceding
$|Z|=r-2$ case. Thus the new maximizing cycle still has one marked
path. Its unchanged unused family recovers the same set $Z$ and the
same unmarked segment $S$, including the order of the core pairs
adjacent to $Z$. Consequently all $\lambda$ representatives
containing the omitted vertex $u$ lie on the original marked run.
They are therefore all marked representatives of the original cycle,
contrary to the choice of $u$, which lies in at most $r-1$ of them.

Hence a trace exists and $0\le \lambda\le n-2r+1$. Equality
forces each trace to alternate across all unmarked pairs, permitting
at most two traces for $\lambda=0$ and one otherwise. At least
$D(n,r)-\lambda$ are required. This excludes $n<2r$ and equality at $n\ge2r$, apart from
$(n,r,\lambda)=(8,3,3)$.
For $\lambda=1$ and $n=2r+1,2r+2$, Lemma~\ref{turn21:unique-mark-boundary}
confines traces to a core path of order at most $2r-3$, with at most
one independent $(r-1)$-set instead of the required $r$.

Only $n=2r+2$, $\lambda=2$, $k\ge3$ remains outside $(8,3)$.
At least $k$ traces lie on at most $2k$ vertices. By
\eqref{eq:released-local-union}, intersecting traces have union of
size at most $k+1$, so distinct ones meet in $k-1$ vertices.
Two disjoint traces would exhaust the vertices, forcing a third
to contain at least $k-1$ from each, impossible. Any three traces
therefore share a vertex and have union of size $k+1$ by
\eqref{eq:released-local-union}. Their three $k$-subsets
cover every pair of the union, giving the $K_r$ excluded by
Lemma~\ref{lem:independent-clique} with $t=r$.
\end{proof}

\begin{proposition}\label{prop:large-complete}
Let $\mathcal H=(H_1,\ldots,H_n)$ be an $r$-graph system on $n$ vertices,
$3\le r<n\le2r+2$, with $\delta_1(H_i)\ge D(n,r)$ for every $i$.
Then $\mathcal H$ has a rainbow Berge Hamilton cycle, except
possibly at $(n,r)=(8,3)$. A counterexample there has a rainbow
Berge seven-cycle and maximum $\lambda\in\{2,3\}$.
\end{proposition}
\begin{proof}
At $(n,r)=(4,3)$, each layer is complete, and distinct triples
$123,234,134,124$ represent the Hamilton core order $(1,2,3,4)$
in distinct colours. Assume $n\ge5$ and no rainbow Berge Hamilton
cycle exists.

For $n\in\{2r+1,2r+2\}$, $t=r$ and
Proposition~\ref{prop:boundary-long} gives an $(n-1)$-cycle.
Lemma~\ref{lem:near-large} forces $\lambda=0$ except possibly
at $(8,3)$, while Lemma~\ref{new:even-distance-two}(ii) excludes
$\lambda=0$ in all cases. Only $(8,3)$ remains, with $\lambda=2,3$.

\emph{Excluding an $(n-1)$-cycle.} Let $n\le2r$.
By Lemma~\ref{lem:near-large}, such a cycle $C$ requires $n=2r$
and $\lambda=0$.
Its unused colour has at least $r$ traces, all maximum independent
sets of $C_{2r-1}$. Exactly $r-1$ such sets contain a given vertex,
having one cyclic gap three and all others two. Hence a trace
avoids each vertex and meets its neighbourhood, so $N_C(U)=V(C)$.

Lemma~\ref{lem:zero-neighbour-load}, with $t=r-1$ and
$\varepsilon=1$, gives $d_R(v)\le r-1$ throughout $C$.
Since every representative avoids the omitted vertex,
$r(n-1)=\sum_{v\in V(C)}d_R(v)\le(r-1)(n-1)$, a contradiction.

\emph{Rainbow Berge spanning paths.} We show that every endpoint $y$ of a
rainbow Berge spanning path with representatives $R$ has
$d_R(y)\le r-2$.
Fix $y$ and apply Lemma~\ref{rt:lem:orbit} to rotate the opposite
endpoint to a vertex $x$ contained in an unused edge of the unused
colour. Write the resulting spanning path as
$P=v_0h_0v_1h_1\cdots h_{n-2}v_{n-1}$, where
$v_0=x$ and $v_{n-1}=y$. Mark the core pair $v_iv_{i+1}$ exactly
when its selected representative $h_i$ contains $y$.

The first core pair $v_0v_1$ cannot be marked: if $y\in h_0$,
then $h_0$ can be assigned to $yv_1$ and closes the path
$v_1h_1v_2\cdots h_{n-2}y$ to a rainbow Berge $(n-1)$-cycle
on $V\setminus\{x\}$, contrary to the preceding result.
Every representative immediately following a marked pair also avoids
$x$. Indeed, if $h_{i+1}$ contained $x$, rotating with $h_{i+1}$ on
$xv_{i+2}$ would make the pair represented by $h_i$ the first pair,
contrary to the preceding observation.

Likewise, an unused edge $f$ at $x$ in the unused colour avoids each
marked head and its successor. If $f$ contained a marked head
$v_{i+1}$, then using $f$ on $xv_{i+1}$ and $h_i$ on $yv_i$ would
close a rainbow Berge Hamilton cycle. If $f$ contained the successor
$v_{i+2}$, using $f$ on $xv_{i+2}$ would rotate the path so that
$h_i$ is the first representative, again contradicting that the
first pair cannot be marked.

Let $T$ be the complement of these forbidden vertices. It contains
the first pair and an unused $r$-edge at $x$, so $|T|\ge r$.
If $|T|=r$, assign $T$ to the first core pair $v_0v_1$ in the
unused colour. Both $T$ and the original representative $h_0$ of that
pair avoid $y$, so the marked pairs and forbidden set remain unchanged.
The exchange releases $h_0$ in its former colour, which is now the
unused colour. The same forbidden-vertex argument therefore gives
$h_0\subseteq T$. Since $h_0$ and $T$ are distinct $r$-sets, this is
impossible. Thus $|T|\ge r+1$.

If $d_R(y)\ge r-1$, then $r+1\le |T|\le n-d_R(y)$ and
$n\le2r$ force $n=2r$, $d_R(y)=r-1$ and $|T|=r+1$. Equality in the count of
forbidden vertices forces the marked pairs to form the terminal
segment beginning at $v_r$, and $T=\{v_0,\ldots,v_r\}$.
Every $v_i$, $0\le i<r$, is attainable as first end with this
terminal segment fixed. For $i>0$, rotate at an unused edge on
$xv_{i+1}$. If none exists, the edge at $x$ is
$T\setminus\{v_{i+1}\}$; putting it on the first pair releases a
different $r$-subset of $T$, which contains $v_{i+1}$ and permits
the rotation in its released colour.

Every representative after the first marked pair must therefore
avoid $v_0,\ldots,v_{r-1}$. All $r-2$ such representatives then equal
$\{v_r,\ldots,v_{2r-1}\}$, contradicting distinctness if $r\ge4$.
For $r=3$, the preceding case gives $n=2r=6$, so $y=v_5$.
The first marked pair is then $v_3v_4$, and its representative
contains $v_3,v_4,v_5$, hence equals $\{v_3,v_4,v_5\}$.
The next representative must avoid $v_0,v_1,v_2$ and contains
$v_4,v_5$, so it is also $\{v_3,v_4,v_5\}$, contrary to the
distinctness of the representatives. This proves the endpoint bound.

\emph{A maximal nonspanning path.} Choose a nonspanning rainbow
Berge path of maximum order $m$, then maximize $\Lambda=d_R(w)$
over such paths and omitted vertices $w$. We have $m\ge3$, since a two-vertex Berge path extends in a second
colour. No rainbow Berge $m$-cycle exists. For $m=n-1$ this was proved
above. Suppose $m<n-1$ and let $C$ be a hypothetical rainbow Berge
$m$-cycle. If one of its representative hyperedges meets
$V\setminus V(C)$, open $C$ at the core pair represented by that
hyperedge and use it to append an exterior vertex. This gives a
nonspanning rainbow Berge path of order $m+1$, contrary to the
maximality of $m$. Thus every representative of $C$ is contained
in $V(C)$.

For each layer $H_i$, its ordinary shadow is connected. Indeed,
every component has at least $r+1$ vertices: if a component had at
most $r$ vertices, each of its vertices would lie in at most one
$r$-edge, contradicting $\delta_1(H_i)\ge r$. Since $n\le2r$,
we have $r+1>n/2$, so there cannot be two components. In an unused
colour $a$, connectedness gives a shadow pair $xy$ with
$x\in V(C)$ and $y\notin V(C)$. By the definition of the ordinary shadow,
some hyperedge $f\in H_a$ contains $x,y$. Since all representatives
of $C$ lie in $V(C)$, this crossing hyperedge is not in $R$ and is
available. Open $C$ at a core pair incident with $x$, and use $f$
to append $y$. This again gives a nonspanning rainbow Berge path
of order $m+1$, a contradiction.

An unused edge at an endpoint $x$ in an unused colour avoids the
opposite end $z$, or closes an $m$-cycle. It also lies in the core
when $m<n-1$, or extends the path.

For $m=n-1$, an unused edge $f$ containing $xw$ instead gives
a rainbow Berge spanning path. The endpoint bound at $w$ then gives
$\Lambda+1\le r-2$. Rotate the original path with $f$ at each
pair entering a vertex of $f\setminus\{x,w\}$. These $r-2$
rotations give paths of order $m$ still omitting $w$. Maximality
of $\Lambda$ forces each discarded representative to contain $w$.
They are distinct, so $\Lambda\ge r-2$, a contradiction.

The path has at least two unused colours; fix two of them,
$a,b$. By Lemma~\ref{rt:lem:orbit}, rotate one endpoint to
$x\in f$ for some $f\in H_a\setminus R$. Then fix $x$ and
rotate the opposite endpoint to $z\in g$ for some
$g\in H_b\setminus R$. Both rotations preserve $R$, its colours and $\Lambda$.

Lemma~\ref{lem:common-endpoint-cuts}(i) gives disjoint sets of
$r-1$ cuts whose heads are determined by $f$ and $r-1$ cuts
whose tails are determined by $g$. Since
$2r-2\le m-1\le n-2\le2r-2,$
we have $n=2r$ and $m=n-1$. The two sets partition the cuts, so either edge determines the
other. Hence $g$ is the unique unused colour-$b$ edge at $z$,
and $d_R(z)\ge r-1$.

The first cut belongs to the set determined by $f$, since
$x\notin g$. Let $h$ be its representative. If $w\in h$,
assign $f$ in colour $a$ to the first pair and use $h$ in its
original colour on $wx$. This gives a rainbow Berge spanning path retaining $R$, with
selected degree at least $r-1$ at $z$, contrary to the endpoint bound.

If $w\notin h$, replace $h$ on the first pair by $f:a$.
Both avoid $w$, so $\Lambda$ is preserved; $g$ remains unused in colour $b$. The cut partition for $g$ and the released $h$ in its
former colour gives $h=f$, since they determine the same head cuts.
This contradicts the original distinctness.
\end{proof}

\section{Small orders with \texorpdfstring{$4\le r<t$}{4 at most r below t}}\label{sec:small-orders}

Before Lemma~\ref{lem:places-count}, we fix the notation used in it.
Here $C$ is a maximizing rainbow Berge $(n-1)$-cycle, $R=E(C)$,
$m=n-1$, and $U$ is its trace union. For each $A\subseteq V(C)$,
retain the orientations $ab$ of core pairs with $a\in U$ and
$b\notin A$. Let $X_A$ be the set of core pairs having a retained
orientation, and let $Z_A$ be the set of retained heads whose two
incident core pairs belong to $X_A$. Put
\[
 T_A=\sum_{a\in A}(d_R(a)-2),
\]
and let $c_A$ count incidences of representatives on $X_A$ with
vertices of $A$ outside their represented pair.

\begin{lemma}\label{lem:places-count}
Let $\mathcal H=(H_1,\ldots,H_n)$ be an $r$-graph system on
$n=2t+1+\varepsilon$ vertices, $t\ge3$ and $\varepsilon\in\{0,1\}$, with
$3\le r<n$ and $\delta_1(H_i)>\binom t{r-1}$ for every $i$.
Assume $\mathcal H$ has no rainbow Berge Hamilton cycle, and let
$C$ be a maximizing rainbow Berge $(n-1)$-cycle. Suppose $r=3$
or $r<t$. For any $A\subseteq V(C)$, form $X=X_A$ and $Z=Z_A$
as above. Assume $d_R(b)\le\binom{t-1}{k-1}$ for every retained
head $b$, whether or not $b\in Z$. Then
\begin{equation}\label{eq:places}
 |Z|+T_A-c_A+\lambda\le(k-1)(m-|X|).
\end{equation}
\end{lemma}
\begin{proof}
Fix a retained orientation $ab$, and release its representative
$h_{ab}$ in its assigned colour $i$. By
Lemma~\ref{lem:protected-marked-cut}, the neighbourhood $N$ at $b$
has size at most $t$. At least $(\binom tk+1)-d_R(b)+1$ edges are available, exceeding
$\binom{t-1}k$ by at least two. Hence $|N|=t$, with two choices
per neighbour, including an alternative for $ab$.
Each pair in $X$ has a retained orientation $ab$ with $a\in U$.
Since $U$ avoids marked ends, the representative currently assigned
to $ab$ avoids $w$. A same-colour exchange using a representative containing $w$ would
increase $d_R(w)$ and contradict maximality. Thus both the original
representative and every allowed replacement avoid $w$, preserving
the marks and the trace union $U$.

Now fix $b\in Z$ and one retained orientation $ab$.
Degree supplies a selected edge $g\in H_i(b)$ leaving $N\cup\{b\}$;
necessarily $g\ne h_{ab}$. Suppose its core pair belongs to $X$,
and choose a retained orientation $uv$ of that pair, with assigned
colour $j$. Here $b$ need not be an endpoint of $uv$.
The hypothesis applies to the retained head $v$, so the preceding
two-choice argument at $v$ supplies an edge $g'\in H_j\setminus R$
containing $u,v$. Replace $g$ by $g'$ on $uv$ in colour $j$.
Both edges avoid $w$, so the cycle remains maximizing and $U$ is
unchanged. In the colour-$i$ neighbourhood released at $b$, this
exchange makes $g$ available and can remove only $g'$. Every old
neighbour in $N$ had at least two available representatives, so
none is lost. Since $g$ leaves $N\cup\{b\}$, the new neighbourhood
has more than $t$ vertices, contrary to
Lemma~\ref{lem:protected-marked-cut}.
Thus the core pair of $g$ lies outside $X$. Since $b\in Z$, both
core pairs incident with $b$ lie in $X$, so $g$ contains $b$ outside
its represented pair.
The representatives outside $X$ also account for $T_A-c_A$
such incidences in $A$ and $\lambda$ at $w$. These vertex sets are
disjoint, and each of the $m-|X|$ representatives has $k-1$ vertices
outside its pair, proving \eqref{eq:places}.
\end{proof}

\begin{proposition}\label{new:three-odd-exceptions}
Let $\mathcal H=(H_1,\ldots,H_n)$ be an $r$-graph system on $n$ vertices,
with $\delta_1(H_i)\ge D(n,r)$ for every $i$.
Then $\mathcal H$ has a rainbow Berge Hamilton cycle if $r=4$ and
$12\le n\le26$, if $r=5$ and $n\in\{13,14,16,17,18\}$,
or if $r=6$ and $n\in\{16,18\}$.
\end{proposition}
\begin{proof}
Suppose otherwise. Proposition~\ref{prop:boundary-long} gives a
maximizing cycle $C$; put $u=|U|$. Lemma~\ref{merge:uniform-zero}
settles $(n,r)\in\{(13,5),(14,5),(16,6),(18,6)\}$ and forces
$\lambda=0$ in the remaining cases unless $r=4$, $t\ge7$.

Apply \eqref{eq:places} with $A=\varnothing$. Its degree hypothesis
holds at every core neighbour $b$ of $a\in U$: if $\lambda=0$,
Lemma~\ref{lem:zero-neighbour-load} gives $d_R(b)\le t$;
otherwise $k=3$, $t\ge7$, and Lemma~\ref{lem:released-exclusions}
gives $d_R(b)\le m-d_\Gamma(a)\le2t-1$.
Both bounds are at most $\binom{t-1}{k-1}$.

The proof of Lemma~\ref{lem:places-count} supplies a core pair
disjoint from $U$. Let $c$ count the core runs of $U$, $s$ the singleton runs,
and $g$ the singleton gaps. Then $|X|=u+c$ and $|Z|=u-s+g$. Equation~\eqref{eq:places}
and the gap lengths give
\[
 ku+(k-1)c-s+g+\lambda\le(k-1)m,\qquad 2c-g\le m-u.\tag{*}
\]

\smallskip\noindent\emph{Case 1: $\lambda=0$.}
Degree gives $u\ge t+1$, so $s<c$. By
Lemma~\ref{new:even-distance-two}(i), traces have mutual core
distances at least three. Write $C_m$ for the ordinary core cycle
and $C_m^2$ for its square, joining vertices at distance at most two.
For $r=4$, traces are independent triples of $C_m^2[U]$.
Delete from $C_m^2[U]$ the $g$ edges joining the endpoints of
gaps of length one. The resulting ordinary graph $F$ is a disjoint
union of squares of paths, with $e(F)=2u-3c+s$.
By inclusion--exclusion, the number of independent triples in $F$ is
\[
 \binom u3-(u-2)e(F)+\sum_{v\in U}\binom{d_F(v)}2-c_3(F).
\]
For a run of length $l$, the last two terms contribute
$B(l)=\max\{0,2l-4,5l-14\}.$
For $l=1,2,3$ these are $0,0,2$. For $l\ge4$, the degrees are
$2,3,4,\ldots,4,3,2$, with $l-2$ consecutive triangles.
Convexity and $B(2)=0$ concentrate all excess run lengths in one
run, bounding the total by $B(u-2c+s+2)$.

Restore the $g$ deleted edges one at a time. Each endpoint then has degree at most two: a run end has at most
two neighbours within its run, and a singleton borders at most one
other restored gap. At least $u-6$ vertices are nonadjacent to both
ends, so each addition destroys at least $u-6$ independent triples.
The number of traces is therefore at most
\[
 \binom u3-(u-2)(2u-3c+s)+B(u-2c+s+2)-g(u-6).
\]

By $(*)$, $t+1\le u\le\lfloor4t/3\rfloor$.
Maximizing over the admissible parameters gives these bounds,
in increasing order of $u$:
\[
\begin{array}{c|l@{\qquad}c|l}
t&\text{upper bounds}&t&\text{upper bounds}\\\hline
5&3&9&63,74,56\\
6&11,4&10&99,106,101\\
7&23,19&11&136,158,158\\
8&41,42&12&193,210,228,220 .
\end{array}
\]
To reproduce the table, put $M=2t+1$ and $1\le c\le M-u$.
Replacing $g$ by its relaxed minimum $g_0=\max\{0,2c+u-M\}$
does not decrease the bound. By $(*)$ and $u\ge2c-s$, it suffices
to consider $s\in[s_0,c-1]$, where
$s_0=\max\{0,3u+2c+g_0-2M,2c-u\}.$
A core pair disjoint from $U$ gives $g<c$. Discard $g_0\ge c$
or $s_0\ge c$; convexity leaves only $s=s_0,c-1$ for each $c$. At $t=12,u=15$, the bounds for $c=1,2,3,4$ are
$165,194,215,228$. Only $228$ can reach $D(n,r)$, and attaining
it requires $s=3$, $g=0$ by $(*)$.
Each singleton run is an isolated vertex in $C_m^2[U]$, whereas
\eqref{eq:uniform-released-degree} bounds its degree in $\Gamma$ by
twelve. Choose a nonneighbour in $\Gamma$ from the other fourteen
vertices of $U$. Their pair lies in at least $u-2-4=9$ independent triples of
$C_m^2[U]$, none of which is a trace. Thus $228$ decreases to $219$. All bounds are below
$D(n,r)=\binom t3+1$.

For $r=5$, a path of order $a$ has $P(a,j)=\binom{a-2j+2}j$
$j$-sets at mutual distance at least three: subtract $2(i-1)$
from the $i$th index. Put $P(a,j)=0$ for $a<3j-2$.
Delete a longest gap from the ordinary core cycle. When $t=7$,
the resulting path gives at most $P(13,4)=\binom74$ traces.
When $t=8$, it gives at most $P(14,4)=\binom84$ traces unless
$m=17$ and every gap has length at most two.
In that exceptional configuration, $g=2c+u-17$, and $(*)$ yields
$5u+4c\le67$. Since $u\ge9$, the number of gaps of length two
is $17-u-c\ge3$.

Choose three such gaps and delete the two bounding a longest
intervening path. The three intervening paths have total order eleven,
so the longest has at least four vertices; the other remaining path
contains the third gap and two nonempty intervening paths, and also
has at least four.
Their orders are $a,13-a$, with $a=4,5,6$ up to interchange.
The bounds $\sum_{j=0}^4P(a,j)P(13-a,4-j)$ are $61,65,66$,
respectively, all below $D(n,r)=71$.

\smallskip\noindent\emph{Case 2: $\lambda>0$.}
Here $r=4$ and $7\le t\le12$. Since $D(n,r)>m$, the trace family
is nonempty. Equation~\eqref{eq:uniform-released-degree} gives
$2\le\Delta(\Gamma)\le t+\varepsilon-1-\lambda$, and hence
$\lambda\le t-2$.
For $u>t$, some core run has length at least two, giving $s<c$
and $2c-s\ge2$. Thus $(*)$ gives $3u\le2m-\lambda-2\le4t-1$,
so $u\le t+3$ for $t\le12$. The same bound holds for $u\le t$.

If $\lambda\ge1+\varepsilon$, then $\Delta(\Gamma)\le t-2$.
Add isolated vertices to reach order $t+3$. The bounded-degree
clique count~\cite[Proposition~2.1]{GLS}, with $\Delta=t-2$ and
$b=4\le\Delta+1$, gives
\[
 c_3(\Gamma)\le\binom{t-1}3+\binom43
 <\binom t3+1-(t-2)\le D(n,r)-\lambda ,
\]
since $\binom{t-1}2>t+1$ for $t\ge7$.
In the remaining case $\varepsilon=\lambda=1$,
Lemma~\ref{lem:independent-clique} gives $K_t\nsubseteq \Gamma$.
Lemma~\ref{new:near-complete-triangle-capacity} with $k=q=3$
gives $c_3(\Gamma)\le\binom t3-1<D(n,r)-1$.
\end{proof}

\section{\texorpdfstring{$3$-graph systems}{3-graph systems}}\label{sec:r3odd}
\begin{proposition}\label{r3:prop:terminal}
Let $\mathcal H=(H_1,\ldots,H_n)$ be a $3$-graph system on $n\ge8$
vertices with $\delta_1(H_i)\ge D(n,3)$ for every $i$.
If it has a rainbow Berge $(n-1)$-cycle, then it has a rainbow
Berge Hamilton cycle.
\end{proposition}
\begin{proof}
Suppose otherwise. Put $m=n-1$ and choose a maximizing rainbow
Berge $m$-cycle $C$, with representatives $R$, omitted vertex $w$
and unused colour $\alpha$. Define the ordinary link graph $L$ on $V(C)$ by
\[
 E(L)=\{ab:\{w,a,b\}\in H_\alpha\setminus R\}.
\]
Put $u=|U|$.
The \emph{third vertex} of a representative triple is the vertex
outside its assigned pair. Each core vertex $v$ has $d_R(v)-2$
such occurrences, totalling $m-\lambda$. We first prove
\begin{equation}\label{r3:eq:basicbounds}
 u\ge t+1,\qquad e(L)\ge D(n,3)-\lambda,\qquad
 \Delta(L)\le t+\varepsilon-1-\lambda,\qquad
 \lambda\le t-2+\varepsilon.
\end{equation}
Degree gives $e(L)\ge D(n,3)-\lambda$. If $L$ is empty,
$\lambda\ge D(n,3)$, excluded at $n=8$ by
Proposition~\ref{prop:large-complete}. For $t\ge4$, at least
$m-2$ pairs would be marked. An unused colour-$\alpha$ triple
avoids $w$ and has at most one marked tail by Lemma~\ref{lem:released-exclusions}. The at most
two unmarked tails allow at most $m-2$ such triples, contradicting
$|H_\alpha\setminus R|\ge nD(n,3)/3-m>m-2$.

For $a\in U$, put $b=a^+$, $s=d_L(a)$ and
$M=m-1-s-\lambda$. Equation~\eqref{eq:released-exclusions}
bounds the released neighbourhood at $b$ by $M$. The $s$ representatives leaving $N_L(a)$ are unmarked and avoid
$b$; among marks, only $bb^+$ can contain $b$. Hence
\[
 D(n,3)\le\binom M2+d_R(b)-1\le\binom{M+1}2+\ind_{bb^+\text{ marked}}.
\]
It follows that $M\ge t-1$. Equality requires $bb^+$ to be
marked and $d_R(b)\ge t+1$. In that case, replace the core
vertex $b$ by $w$, using an unused triple containing $aw$ and
the marked representative on $wb^+$. Global maximality then gives
$\lambda\ge d_R(b)-1\ge t$. The inequality
$1\le s\le m-t-\lambda$ forces $m=2t+1$, $\lambda=t$ and
$d_L(a)=1$ for every $a\in U$.

Apply equality in both orientations. The core neighbours of an edge
of $L$ include three vertices of selected degree at least $t+1$,
requiring $3(t-1)>m-\lambda$ third-vertex occurrences. Thus
$M\ge t$, proving the maximum-degree bound.

If $u\le t$, then $\lambda=0$ is excluded by $e(L)>\binom t2$.
Otherwise,
$t(t-1)+2-2\lambda\le2e(L)\le t(t-\lambda),$
forcing $\lambda=1$ and $L[U]=K_t$. This contradicts the
maximum-degree bound in odd order and Lemma~\ref{lem:independent-clique}
in even order. Thus $u\ge t+1$. Marked ends give $u\le m-\lambda-1$ if
$\lambda>0$, and $e(L)>0$ gives $\lambda\le t-2+\varepsilon$.
This proves \eqref{r3:eq:basicbounds}.

At $n=8$, Proposition~\ref{prop:large-complete} and
\eqref{r3:eq:basicbounds} give $\lambda=2$, $\Delta(L)\le1$ and
$e(L)\ge2$. Marked ends leave at most four active vertices, forcing
$L[U]=2K_2$, excluded by Lemma~\ref{lem:released-exclusions}.
Henceforth $t\ge4$.

Choose $W\subseteq V(C)$ containing every
$b\in N_C(U)$ with $d_R(b)\ge t$. Lemma~\ref{lem:places-count}
applies with $A=W$, since every retained head has degree
at most $t-1$. In its notation, put $f=|X|+|Z|$. Then
\begin{equation}\label{r3:eq:charge}
 f+T_W-c_W+\lambda\le m .
\end{equation}
Each retained orientation $ab$ has at least $t-d_R(b)$ unused
third-vertex alternatives by the proof of \eqref{eq:places}.
Before excluding heads in $W$, let $\xi$ count the core runs of $U$,
$\nu$ the nonsingleton runs, and $g$ the singleton gaps.
Then $|X|=u+\xi$, $|Z|=u-\xi+\nu+g$, so the initial value is
$F=2u+\nu+g$. If $U=V(C)$, put $F=2m$.

\smallskip\noindent\emph{Case 1: $\lambda=0$.}
Take $W=\{h\in N_C(U):d_R(h)\ge t\}$ and put $q=|W|$.
Lemma~\ref{lem:zero-neighbour-load} gives $q=0$ in odd order,
but then $f\ge2u>m$ contradicts \eqref{r3:eq:charge}.
Thus $m=2t+1$, $q>0$, and the same lemma gives $d_R(h)=t$
with exactly one active core neighbour for each $h\in W$.
In particular, $U\ne V(C)$.

Keeping the core order, colours and unused link $L$, minimize $q$.
A representative counted by $c_W$ has third vertex in $W$ and pair
in $X$. Every same-colour unused replacement has third vertex
$v\in N_C(U)$ with $d_R(v)=t-1$, or it decreases $q$ or violates
Lemma~\ref{lem:zero-neighbour-load}.

The alternative third vertices are distinct. Replacements sharing
$v$ could be simultaneous, since their core pairs differ, and
would give $d_R(v)=t+1$ in $N_C(U)$.
These vertices contribute $c_W(t-3)$ third-vertex occurrences,
and $W$ contributes $q(t-2)\ge c_W$. Hence $(t-2)c_W\le m$.

Exclude heads in $W$ one at a time. For $h$ with unique active
neighbour $x$, only $hx$ can leave $X$. If it stays, only $h$
can leave $Z$. If it leaves, at most one of $x,h$ was in $Z$: if $h$ is active,
$x$ was already excluded; otherwise the other pair at $h$ lies outside $X$. Thus $f\ge F-2q$.

Put $\ell=F-m$. Since $T_W=q(t-2)$,
\eqref{r3:eq:charge} gives
$\ell+q(t-4)\le c_W\le q(t-2).$
Hence $q\ge\ell/2$, $c_W\ge\ell(t-2)/2$, and
$(t-2)^2\ell\le2m$.

The $2u-2\xi-\nu+g$ pairs of $U$ at core distance at most two
are absent from $L$ by Lemmas~\ref{lem:released-exclusions}
and~\ref{new:even-distance-two}(i). Since $2\xi-g\le m-u$,
$e(L)\ge D(n,3)$ gives
$\nu\ge3u-m-\binom u2+D(n,3).$
If $u=t+1$, then $\nu\ge3$; if $u\ge t+2$, then $\nu\ge1$.
Thus $\ell\ge4$, and $(t-2)^2\ell\le4t+2$ forces $t=4$, $\ell=4$.

The case $u=5$ violates $2\nu\le u$, so $u=6$, $\nu=1$, $g=0$.
The active vertices form one run, with only $\binom42=6$ pairs at
core distance at least three. Thus $e(L)\le6<D(n,3)=7$.

\smallskip\noindent\emph{Case 2: $\lambda>0$ and $t\ge5$.}
Set $W=\{v\in V(C):d_R(v)\ge t\}$, put $q=|W|$, and let
$p=|\{v\in V(C):d_R(v)=t-1\}|$.
There is a nonsingleton run since $u>m/2$, and a nonsingleton gap
containing marked ends. Each excluded head removes at most two
pairs, so $|X|\ge u+\xi-2q$.

Starting after a nonsingleton gap, select all but the first vertex
of each run, plus the first vertex of the first nonsingleton run.
Of these $u-\xi+1$ vertices, an internal vertex leaves $Z$ only if
it lies in $W$; an endpoint leaves only if it or its inactive
neighbour lies in $W$. These inactive neighbours are distinct:
the extra first vertex follows the initial gap or a singleton run.
Thus $|Z|\ge u-\xi+1-q$.

All marks are among the $m-u-\xi$ pairs avoiding $U$; when
$\lambda=1$, Lemma~\ref{turn21:unique-mark-boundary} gives at
least three such pairs. Thus $m-u-\xi+\lambda\ge4$.
Each vertex of $Z$ needs one third-vertex occurrence. A vertex of
degree $t-1$ needs $t-3$, at most one counted for $Z$. Therefore
\[
 T_W+(t-4)p\le m-\lambda-|Z|\le 2t-5+q .
\]
Since $T_W\ge q(t-2)$ and $t\ge5$, we obtain $q\le2$.
Moreover, $f\ge2u+1-3q$, so \eqref{r3:eq:charge} implies
$c_W\ge2u+1-m+\lambda+q(t-5)>0.$
In particular, $q\ge1$.

Among maximizing cycles, minimize first
\[
 \sum_{v\ne w}\max\{0,d_R(v)-t+1\},
\]
then $q$. A representative counted by $c_W$ has third vertex
$h\in W$ and a retained head $b$ of degree $d=d_R(b)<t$.
Same-colour replacement avoids $w$ and preserves $\lambda$.
A third vertex of degree at most $t-2$ decreases the sum; if
$d_R(h)=t$, one of degree at least $t$ decreases $q$ with the sum fixed. Hence all $t-d$ alternatives have degree at
least $t-1$, exactly $t-1$ when $d_R(h)=t$.

Suppose $q=2$, and put $s=T_W-2(t-2)$, the total excess of
$d_R(v)$ above $t$ on $W$. The incidence bound gives
$p+s\le(t-4)p+s\le1$.
An alternative has degree $t-1$ or is the other vertex of $W$,
the latter requiring $s\ge1$. Thus $t-d\le p+s\le1$, so $d=t-1$.
But $b$ itself counts towards $p$ and cannot be its pair's third
vertex, leaving no alternative.

It remains to consider $q=1$, say $W=\{h\}$.
Since $c_W\le T_W$, the bounds $2u-2\le f\le m-\lambda$
force $\lambda=\varepsilon=1$, $u=t+1$ and $c_W=T_W$.
Equality in $|X|\ge u+\xi-2$ requires both core pairs incident
with $h$ to be removed from $X$. Thus $h$ is inactive, while
$x=h^-$ and $z=h^+$ are active.
At most $t$ pairs of $U$ are absent from $L$, including $u-\xi$
core pairs, leaving at most $\xi-1$ missing noncore pairs.

A representative counted by $c_W$ on a pair leaving
$a\in U\setminus\{x\}$ forces $ax\notin E(L)$ by
\eqref{eq:released-exclusions}; one entering
$b\in U\setminus\{z\}$ forces $bz\notin E(L)$ by reversal.
Each missing noncore pair other than $xz$ accounts for at most one
counted representative. Missing core pairs at $x,z$ account only
for $x^-x,zz^+$; both also require $xz\notin E(L)$.
The pair $xz$ accounts for at most these two, giving $c_W\le \xi$.
By Lemma~\ref{turn21:unique-mark-boundary}, the $t$ inactive
vertices include four consecutive core vertices, so $\xi\le t-3$.
This contradicts $T_W=c_W\ge t-2$.

\smallskip\noindent\emph{Case 3: $\lambda>0$ and $t=4$.}
If $\lambda\ge2+\varepsilon$, then
\eqref{r3:eq:basicbounds} gives $\Delta(L)\le1$ and
$\lambda\le3$. Consequently,
$2(7-\lambda)\le u\le m-\lambda-1$, which is impossible.
At $n=9$, only $\lambda=1$ remains, but
Lemma~\ref{turn21:unique-mark-boundary} leaves at most four active
vertices. Hence $n=10$, $\lambda\in\{1,2\}$.

Write the core cycle as
\[
 C=v_0h_0v_1h_1\cdots v_8h_8v_0,
\]
with indices modulo nine. Suppose first that $\lambda=1$, with
marked pair $v_0v_1$.
Lemma~\ref{turn21:unique-mark-boundary}, trace independence and
$e(L)\ge6$ force $U=\{v_3,v_4,v_5,v_6,v_7\}$ and make $L[U]$
the complement of the ordinary path with vertex order
$(v_3,v_4,v_5,v_6,v_7)$.
Applying \eqref{eq:released-exclusions} in both orientations forces
the third vertex of the representative on $v_3v_4$ to be $v_2$,
and that of the representative on $v_6v_7$ to be $v_8$. It also
permits $v_4$ or $v_6$ as a third vertex only on the pairs
$v_1v_2$ and $v_8v_0$.

A core pair oriented from a vertex of $U$ to $b$ with
$d_R(b)\le3$ has at least $5-d_R(b)\ge2$ third-vertex choices
by the proof of \eqref{eq:places}.
The unique choices on $v_3v_4,v_6v_7$ force
$d_R(v_4),d_R(v_6)\ge4$.
Their four required third-vertex occurrences cannot fit on
$v_1v_2,v_8v_0$.

Now let $\lambda=2$. On any maximizing cycle with disjoint marks,
the four marked ends are inactive. The link on the other five
vertices has at least five edges and maximum degree two, so it is
a five-cycle. Degree $D(n,3)=7$ then forces both marked triples
into the unused layer.

Choose an orientation and label a marked pair $v_0v_1$ so that
$v_8$ is active; this is possible because $u\ge5$.
Replace core vertex $v_0$ by $w$, using an unused triple containing
$v_8w$ and the marked representative on $wv_1$.
Trace independence ensures that the inserted triple avoids $v_0$.
The released colour $i$ of $v_8v_0$ is now unused. Its link union
at $v_0$ equals the released neighbourhood $T$, with
$|T|\le4$ by Lemma~\ref{lem:protected-marked-cut}.

A second representative containing $v_0$ would make the new cycle
maximizing, contradicting \eqref{r3:eq:basicbounds} since its link
union has size at most four. Thus only the representative of
$wv_1$ contains $v_0$. Degree seven forces $|T|=4$ and all six
pairs of $T$ into the colour-$i$ link at $v_0$.

If a core pair $xy$ lies in $T$, replacing its representative by
$\{v_0,x,y\}:i$ gives disjoint marks $xy,wv_1$ at the maximum
$\lambda=2$. By the disjoint-mark observation,
$\{v_0,x,y\}$ also belongs to the former layer of $xy$.
Assign it that colour instead, preserving $\lambda=2$ and leaving
$i$ unused with link union contained in $T$. This contradicts
\eqref{r3:eq:basicbounds}.

Hence $T$ is the unique independent four-set of the ordinary core
path $(v_2,v_3,\ldots,v_8)$:
$\{v_2,v_4,v_6,v_8\}$.
Equation~\eqref{eq:released-exclusions} makes $T$ avoid the head of
the other marked pair, so that pair is not $v_1v_2$; it is not
$v_8v_0$ since $v_8$ is active. Thus the marked pairs are disjoint.
If $v_2$ were active, reversing the orientation at $v_0v_1$ would
give the same four-set avoiding the tail of the other marked pair.
Both endpoints would then lie in $\{v_3,v_5,v_7\}$, which contains
no core pair. Hence $v_2$ is inactive. The four inactive vertices
are exactly the marked endpoints, so the marked pairs are
$v_0v_1$ and $v_2v_3$.

In the five-cycle $L[\{v_4,v_5,v_6,v_7,v_8\}]$, trace independence
and degree two give $N_L(v_6)=\{v_4,v_8\}$.
Equation~\eqref{eq:released-exclusions} gives
$N_L(v_8)^+\cap T=\varnothing$, hence
$N_L(v_8)=\{v_4,v_6\}$.
These neighbourhoods force a triangle on $\{v_4,v_6,v_8\}$,
a contradiction.
\end{proof}

\begin{proof}[Proof of Theorem~\ref{thm:main}]
For $n\le2r+2$, Proposition~\ref{prop:large-complete} proves the
claim except at $(8,3)$, where its seven-cycle extends by
Proposition~\ref{r3:prop:terminal}.

Let $n>2r+2$ and put $k=r-1$, so $r<t$. Suppose no rainbow
Berge Hamilton cycle exists. Proposition~\ref{prop:boundary-long}
gives an $(n-1)$-cycle, contradicting Proposition~\ref{r3:prop:terminal}
if $r=3$. Thus $r\ge4$. Lemma~\ref{merge:uniform-zero} excludes
$2k\ge t+1+\varepsilon$, so $2k\le t+\varepsilon$.
If $r=4$, $t\le12$, then $12\le n\le26$.
If $r\ge5$, $t\le8$, then $8\le2k\le t+\varepsilon\le9$
forces $r=5$, $n\in\{16,17,18\}$.
Proposition~\ref{new:three-odd-exceptions} covers both ranges.
Otherwise \eqref{eq:shadow-range} holds, and
Proposition~\ref{prop:shadow-conclusion} completes the proof.
\end{proof}

\section*{Declaration of generative AI use}
OpenAI ChatGPT assisted with proof exploration, organization, language
editing, and checks of arguments and references. The authors are
responsible for verifying all arguments and references and for the
manuscript's content.

\section*{Acknowledgements}
The research is supported by the National Natural Science Foundation of China
(No. 12171088).

\end{document}